\documentclass[10pt]{amsart}
\usepackage[latin1]{inputenc}
\usepackage{amssymb}
\usepackage{amsmath}
\begin{document}
\setlength{\parindent}{1.2em}
\def\COMMENT#1{}
\def\TASK#1{}
\def\noproof{{\unskip\nobreak\hfill\penalty50\hskip2em\hbox{}\nobreak\hfill%
       $\square$\parfillskip=0pt\finalhyphendemerits=0\par}\goodbreak}
\def\endproof{\noproof\bigskip}
\newdimen\margin   
\def\textno#1&#2\par{%
   \margin=\hsize
   \advance\margin by -4\parindent
          \setbox1=\hbox{\sl#1}%
   \ifdim\wd1 < \margin
      $$\box1\eqno#2$$%
   \else
      \bigbreak
      \hbox to \hsize{\indent$\vcenter{\advance\hsize by -3\parindent
      \it\noindent#1}\hfil#2$}%
      \bigbreak
   \fi}
\def\proof{\removelastskip\penalty55\medskip\noindent{\bf Proof. }}
\newtheorem{firstthm}{Proposition}
\newtheorem{thm}[firstthm]{Theorem}
\newtheorem{prop}[firstthm]{Proposition}
\newtheorem{lemma}[firstthm]{Lemma}
\newtheorem{cor}[firstthm]{Corollary}
\newtheorem{problem}[firstthm]{Problem}
\newtheorem{defin}[firstthm]{Definition}
\newtheorem{conj}[firstthm]{Conjecture}
\def\eps{{\varepsilon}}
\def\N{\mathbb{N}}
\def\R{\mathbb{R}}
\def\B{\mathcal{B}}
\def\K{\mathcal{K}}
\def\S{\mathcal{S}}
\def\D{\mathcal{D}}
\def\C{\mathcal{C}}
\def\Cbar{\overline{C}}
\def\Cjbar{\overline{C_j}}
\def\Cprimejbar{\overline{C'_j}}
\def\cH{\mathcal{H}}
\newcommand{\ex}{\mathbb{E}}
\newcommand{\pr}{\mathbb{P}}
\newcommand{\hcf}{{\rm hcf}}

\title{The minimum degree threshold for perfect graph
packings}
\author{Daniela K\"uhn \and Deryk Osthus}
\date{}
\maketitle
\vspace{-.8cm}
\begin{abstract} \noindent
Let $H$ be any graph.
We determine up to an additive constant the minimum degree of a graph~$G$
which ensures that $G$ has a perfect $H$-packing (also called an $H$-factor). More precisely,
let $\delta(H,n)$ denote the smallest integer $k$ such that every graph $G$ 
whose order $n$ is divisible by~$|H|$ and with $\delta(G)\ge k$ contains a
perfect $H$-packing. We show that
$$\delta(H,n) = \left(1-\frac{1}{\chi^*(H)} \right)n+O(1).$$
The value of $\chi^*(H)$ depends on the relative sizes of the colour classes in the 
optimal colourings of $H$ and satisfies $\chi(H)-1<\chi^*(H)\le \chi(H)$.
\end{abstract}

\section{Introduction}\label{intro}
\subsection{Background}
Given two graphs $H$ and
$G$, an \emph{$H$-packing in $G$} is a collection of vertex-disjoint copies
of $H$ in $G$. $H$-packings are natural generalizations of graph matchings
(which correspond to the case when $H$ consists of a single edge). 
An $H$-packing in $G$ is called \emph{perfect} if it
covers all vertices of $G$. In this case, we also say that $G$ contains
an \emph{$H$-factor} or a \emph{perfect $H$-matching}.
If $H$ has a component which contains at least~3 vertices then the question whether 
$G$ has a perfect $H$-packing is difficult from both a structural and algorithmic
point of view: Tutte's theorem characterizes those graphs which have a perfect 
$H$-packing if $H$ is an edge but for other graphs~$H$ no such characterization exists.
Moreover, Hell and Kirkpatrick~\cite{HKsiam}
showed that the decision problem whether a graph $G$ has a perfect $H$-packing
is NP-complete if and only if $H$ has a component which contains at least~3 vertices.
They were motivated by questions arising in timetabling (see~\cite{HK78}).

This leads to the search for simple sufficient conditions which ensure the existence of 
a perfect $H$-packing. A fundamental result of this kind is the theorem of Hajnal and 
Szemer\'edi~\cite{HSz} which states that every graph~$G$ whose order~$n$ is divisible
by~$r$ and whose minimum degree is at least
$(1-1/r)n$ contains a perfect $K_r$-packing. The minimum degree condition is easily 
seen to be best possible. (The case when $r=3$
was proved earlier by Corr\'adi and Hajnal~\cite{CH}.) 
The following result is a generalization of this to arbitrary graphs $H$.

\begin{thm}\label{KSS}{\bf [Koml\'os, S\'ark\"ozy and Szemer\'edi~\cite{KSSz01}]}
For every graph $H$ 
there exists a constant $C=C(H)$ such that every
graph $G$ whose order $n$ is divisible by $|H|$ and whose
minimum degree is at least $(1-1/\chi(H))n+C$ contains a perfect $H$-packing.
\end{thm}

This confirmed a conjecture of Alon and Yuster~\cite{AY96}, who had obtained the
above result with an additional error term of~$\eps n$ in the minimum degree condition.
As observed in~\cite{AY96},
there are graphs $H$ for which the above constant~$C$ cannot be omitted completely.
Thus one might think that this settles the question of which minimum degree 
guarantees a perfect $H$-packing.

However, there are graphs $H$ for which the bound on the minimum 
degree can be improved significantly: Kawarabayashi~\cite{KK} conjectured that
if $H=K_\ell^-$ (i.e.~a complete graph with one edge removed)
and $\ell \ge 4$ then one can replace the chromatic number with the
critical chromatic number in Theorem~\ref{KSS} and take~$C=0$.
He~\cite{KK} proved the case $\ell=4$ and together with Cooley, we proved
the general case for all graphs whose order $n$ is sufficiently large~\cite{KOKlminus}.
Here the \emph{critical chromatic number} $\chi_{cr}(H)$ of a graph $H$
is defined as $(\chi(H)-1)|H|/(|H|-\sigma(H))$, where $\sigma(H)$ denotes
the minimum size of the smallest colour class in a colouring
of $H$ with $\chi(H)$ colours. 
Note that $\chi_{cr}(H)$ always satisfies  $\chi(H)-1 < \chi_{cr}(H) \le \chi(H)$  
and equals $\chi(H)$ if and only if for every colouring of $H$ with $\chi(H)$
colours all the colour classes have equal size.

The critical chromatic number was introduced by Koml\'os~\cite{JKtiling}.
He (and independently Alon and Fischer~\cite{AF99}) observed that for \emph{any}
graph~$H$ it gives a lower bound on the minimum degree that guarantees
a perfect $H$-packing.%
     \COMMENT{This is true for \emph{every} integer $n$ that is divisible by $|H|$.
Indeed, let $\ell:=\chi(H)$ and $\sigma=\sigma(H)$.
Let $k\in\mathbb{N}$ and let $G$ be the complete $\ell$-partite graph of
order $k|H|$ whose smallest vertex class has size $\frac{\sigma}{|H|}k|H|-1=\sigma k-1$
and whose other vertex class sizes are as equal as possible.
Then $G$ doesn't contain a perfect $H$-packing and
$$\delta(G)=k|H|-\text{largest vx class}=
k|H|-\lceil \frac{k|H|-\sigma k+1}{\ell-1}\rceil.$$
But $|H|-\sigma=|H|(1-\xi/(\ell-1+\xi))$ since
$(\ell-1)\sigma=z_1=\xi z=\xi(|H|-\sigma)$. So
\begin{align*}
\delta(G) =k|H|-\left\lceil
\frac{k|H|\left(\frac{\ell-1}{\ell-1+\xi}\right)+1}{\ell-1}\right\rceil
\ge k|H|-\frac{k|H|\left(\frac{\ell-1}{\ell-1+\xi}\right)+1}{\ell-1}-\frac{\ell-2}{\ell-1}
=k|H|\left(1-\frac{1}{\ell-1+\xi}\right)-1.
\end{align*}
(To see the inequality use that $k|H|\frac{\ell-1}{\ell-1+\xi}$ is an
integer.)}

\begin{prop}\label{propKomlos}
For every graph $H$ and every integer $n$ that is divisible by $|H|$ there exists
a graph $G$ of order $n$ and minimum degree $\lceil(1-1/\chi_{cr}(H))n\rceil-1$
which does not contain a perfect $H$-packing.
\end{prop} 

Koml\'os also showed that the critical chromatic number is the parameter which governs the 
existence of \emph{almost} perfect packings in graphs of large minimum degree.

\begin{thm}\label{thmKomlos}{\bf [Koml\'os~\cite{JKtiling}]}
For every graph $H$ and every $\gamma>0$ there exists an integer
$n_1=n_1(\gamma,H)$ such that every graph $G$ of order $n\ge n_1$
and minimum degree at least $(1-1/\chi_{cr}(H))n$ contains an $H$-packing
which covers all but at most $\gamma n$ vertices of~$G$.
\end{thm}

Confirming a conjecture of Koml\'os~\cite{JKtiling}, Shokoufandeh and Zhao~\cite{SZ} proved
that the number of uncovered vertices can be reduced to a constant depending only on~$H$.

\subsection{Main result}
Our main result is that for any graph~$H$, either its critical chromatic number or its
chromatic number is the relevant parameter which governs the existence of perfect 
packings in graphs of large minimum degree. The exact classification depends on a
parameter which we call the highest common factor of $H$ and which is defined as follows.

We say that a colouring of $H$ is \emph{optimal} if it uses exactly $\chi(H)=:\ell$ colours. 
Given an optimal colouring $c$, let $x_1\le x_2\le \dots\le x_{\ell}$
denote the sizes of the colour classes of~$c$. Put
$\D(c):= \{x_{i+1}-x_i\,|\, i=1,\dots,\ell-1\}.$
Let $\D(H)$ denote the union of all the sets $\D(c)$ taken over all optimal colouring $c$.
We denote by $\hcf_\chi(H)$ the highest common factor of all integers in $\D(H)$.
(If $\D(H)=\{0\}$ we set $\hcf_\chi(H):=\infty$.)
We write $\hcf_c(H)$ for the highest common factor of all the orders of components
of~$H$. If $\chi(H)\neq 2$ we say that $\hcf(H)=1$ if $\hcf_\chi(H)=1$.
If $\chi(H)=2$ then we say that $\hcf(H)=1$ if both $\hcf_c(H)=1$ and $\hcf_\chi(H)\le 2$.
The following table gives some examples:
$$
\begin{array}{l|l|l|l|l|l}
H & \chi(H) & \chi_{cr}(H) & \hcf_\chi(H) & \hcf_c(H) & \hcf(H) \\
\hline
C_{2k+1} \ (k\ge 2) & 3 & 2+1/k & 1 & - & 1 \\
C_{2k}  & 2 & 2   & \infty & 2k & \neq 1 \\
K_{1,2}\cup C_6 & 2 & 9/5   & 1 & 3   & \neq 1 \\
K_{1,4}\cup C_4 & 2 & 3/2   & 3 & 1 & \neq 1 \\
K_{1,2}\cup K_{1,4} & 2  & 4/3   & 2 & 1 & 1 
\end{array}
$$
Note that if all the optimal colourings of $H$ have the property that
all colour classes have equal
size, then $\D(H)=\{0\}$ and so $\hcf(H) \neq 1$ in this case.
So if $\chi_{cr}(H)=\chi(H)$, then $\hcf(H) \neq 1$.
Moreover, it is easy to see that
there are graphs $H$ with
$\hcf_\chi(H)=1$ but such that for all optimal colourings $c$ of $H$ the highest
common factor of all integers in $\D(c)$ is strictly bigger than one (for example,
take~$H$ to be the graph obtained from~$K_{1,4,6}$ by adding a new vertex and joining
it to all the vertices in the vertex class of size~4).%
     \COMMENT{Indeed, take for $H$ the graph obtained from $K_{1,4,6}$ by
adding one vertex joined to the vertex class of size~4. Then the optimal
colourings of $H$ have vertex classes of size $2,4,6$ and $1,4,7$ respectively.
The hcf of the 1st colouring is 2 and the hcf of the 2nd one is 3.}
Thus for such graphs~$H$ we do need to consider all optimal colourings of $H$. 

As indicated above, our main result is that in Theorem~\ref{KSS}
one can replace the chromatic number by the critical chromatic number if $\hcf(H)=1$. 

\begin{thm}\label{thmmain}
Suppose that $H$ is a graph with
$\hcf(H)=1$. Then there exists
a constant $C=C(H)$ such that every graph $G$ whose order $n$ is divisible by~$|H|$
and whose minimum degree is at least $(1-1/\chi_{cr}(H))n+C$
contains a perfect $H$-packing.
\end{thm}

Note that Proposition~\ref{propKomlos} shows the result is best possible up to the
value of the constant~$C$. A simple modification of the examples in~\cite{AF99,JKtiling}
shows that there are graphs~$H$ for which the constant $C$ cannot be omitted
entirely.%
      \COMMENT{Indeed, let $H=K_{s,s,s-1}$ where $s\ge 6$. So $\hcf(H)=1$.
Let $G$ be complete tripartite graph with smallest vertex class size
$\sigma k-1$ and other vertex class sizes as equal as possible (i.e. the graph from the 
first comment). Also add an $(s-1)$-factor of large girth into each of the vertex classes.
We claim that $G$ doesn't contain a perfect $H$-packing. Indeed, suppose that $G$
has a perfect $H$-packing. Let $A$ be one of the large vertex classes of~$G$. Then (wlog)
some copy of $H$ in the $H$-packing must have at least $s+1$ vertices in~$A$.
Let $c_i$ denote the size of the intersection of the $i$th colour class of $H$ with~$A$.
Wlog $c_1\le c_2\le c_3$. Note that $c_3\le s-1$ as $G[A]$ is $(s-1)$-regular.
Thus either $c_1,c_2\ge 1$ or $c_2\ge 2$. In both cases we must have a 4-cycle in $G[A]$,
a contradiction.}
Moreover, it turns out that Theorem~\ref{KSS} is already best possible up to
the value of the constant~$C$ if $\hcf(H)\neq 1$ (see Propositions~\ref{bestposs1}
and~\ref{bestposs2} for the details).
In~\cite{KOSODA} we sketched a simpler argument that yields a weaker result than
Theorem~\ref{thmmain}: there $H$ had to be either connected or non-bipartite and
we needed an additional error term of~$\eps n$ in the minimum degree condition.

If we combine Theorems~\ref{KSS} and~\ref{thmmain} together with
Propositions~\ref{propKomlos},~\ref{bestposs1} and~\ref{bestposs2}
we obtain the statement indicated in the abstract. Let 
$$\chi^*(H):=
\begin{cases}
\chi_{cr}(H) &\text{ if $\hcf (H)=1$};\\
\chi(H) &\text{ otherwise}.
\end{cases}
$$
Also let $\delta(H,n)$ denote the smallest integer $k$ such that every graph $G$ 
whose order $n$ is divisible by~$|H|$ and with $\delta(G)\ge k$ contains a perfect $H$-packing.

\begin{thm}\label{thmmaingeneral}
For every graph $H$ there exists a constant $C=C(H)$ such that
$$\left( 1-\frac{1}{\chi^*(H)} \right)n-1\le \delta(H,n)
\le \left(1-\frac{1}{\chi^*(H)} \right)n+C.$$
\end{thm}

\noindent (The $-1$ on the left hand side can be omitted if $\hcf(H)= 1$ or
if $\chi(H)\ge 3$.)
Thus for perfect packings in graphs of large minimum degree, 
the parameter $\chi^*(H)$ is the relevant parameter, whereas for 
almost perfect packings it is $\chi_{cr}(H)$. Note that while the definition of the parameter
$\chi^*$ is somewhat complicated, the form of Theorem~\ref{thmmaingeneral} is exactly analogous
to that of the Erd\H{o}s-Stone theorem (see
e.g.~\cite[Thm~7.1.2]{Diestel} or~\cite[Ch.~IV, Thm.~20]{BGraphTh}), which implies that 
\begin{equation}\label{eqErdosStone}
ex(H,n)=\left(1-\frac{1}{\chi(H)-1} +o(1) \right)n,
\end{equation}
where $ex(H,n)$ denotes the smallest number~$k$ such that every graph $G$ of
order $n$ and average degree~$>k$ contains a copy of~$H$. 

\subsection{Open problems}
Our constant~$C$ appearing in Theorems~\ref{thmmain} and~\ref{thmmaingeneral} is rather large
since it is related to the number of partition classes (clusters) obtained by the Regularity lemma.
It would be interesting to know whether one can take e.g.~$C= |H|$. Another open problem
is to characterize all those graphs~$H$ for which $\delta(H,n)=\lceil(1-1/\chi^*(H))n\rceil$.
This is known to be
the case e.g.~for complete graphs~\cite{HSz} and, if $n$ is large, for cycles~\cite{Abbasi}
and for the case when $H=K_\ell^-$~\cite{KOKlminus}. Further observations on this problem
can be found in~\cite{KOKlminus} and~\cite{MPhil}.

\subsection{Algorithmic aspects}
Kann~\cite{Kann94} showed that the optimization problem of finding a maximum 
$H$-packing is APX-complete if $H$ is connected and $|H| \ge 3$ 
(i.e.~it is not possible to approximate the optimum solution
within an arbitrary factor unless P=NP). For such $H$ and any $\gamma>0$, 
we gave a polynomial 
time algorithm in~\cite{KOSODA} which finds a perfect $H$-packing if 
$\delta(G) \ge (1-1/\chi^*(H) +\gamma )n$.
Also note that Theorem~\ref{thmmain} immediately implies that the decision problem whether
a graph $G$ has a perfect $H$-packing is trivially solvable in polynomial time 
in this case. On the other hand, in~\cite{KOSODA} we showed that for many graphs $H$,
the problem becomes NP-complete when the input graphs are all those graphs $G$ with
minimum degree at least $(1-1/\chi^*(H)-\gamma)|G|$, where $\gamma>0$ is arbitrary.
We were able to show this if $H$ is complete or $H$ is a complete $\ell$-partite graph
where all colour classes contain at least two vertices. It would certainly be interesting
to know whether this extends to all graphs $H$ which have a component with at least three
vertices.

\subsection{Organization of the paper} 
In the next section we introduce some basic definitions and then describe the extremal examples
which show that our main result is best possible.
In Section~\ref{sec:RLandBL} we then state the Regularity lemma of Szemer\'edi
and the Blow-up lemma of Koml\'os, S\'ark\"ozy and Szemer\'edi.
In Section~\ref{sec:overview} we give a rough outline of the structure of the proof
and state some of the main lemmas. In Section~\ref{sec:nonextremal} we consider the case
where $G$ is not similar to
an extremal graph. In Section~\ref{sec:complete} we investigate perfect $H$-packings in
complete $\ell$-partite graphs (these results are needed when we apply the Blow-up lemma
in Sections~\ref{sec:nonextremal} and~\ref{sec:extremal}).
In Section~\ref{sec:extremal} we consider the case where $G$ is similar to an extremal graph.
Finally, we combine the results of the previous sections in Section~\ref{sec:thmproof}
to prove Theorem~\ref{thmmain}. Our proof of Theorem~\ref{thmmain} uses ideas from~\cite{KSSz01}.


\section{Notation and extremal examples}\label{sec:tools}

Throughout this paper we omit floors and ceilings whenever
this does not affect the argument.
We write $e(G)$ for the number of edges of a graph $G$, $|G|$ for its
order, $\delta(G)$ for its minimum degree, $\Delta(G)$ for its maximum
degree, $\chi(G)$ for its chromatic number and $\chi_{cr}(G)$ for its
critical chromatic number as defined in Section~\ref{intro}.
We denote the degree of a vertex $x\in G$ by $d_G(x)$ and its neighbourhood
by $N_G(x)$. Given a set $A\subseteq V(G)$, we write $e(A)$ for the 
number of edges in~$A$ and define the density of~$A$ by $d(A):=e(A)/\binom{|A|}{2}$.

Given disjoint $A,B\subseteq V(G)$,
an \emph{$A$--$B$ edge} is an edge of $G$ with one endvertex in $A$ and the
other in $B$; the number of these edges is denoted by $e_G(A,B)$ or $e(A,B)$
if this is unambiguous. We write $(A,B)_G$ for the bipartite subgraph of $G$
whose vertex classes are $A$ and $B$ and whose edges are all $A$--$B$ edges
in~$G$. More generally, we write $(A,B)$ for a bipartite graph with vertex
classes $A$ and~$B$.

The following two propositions together show that if $\hcf(H)\neq 1$ then
Theorem~\ref{KSS} is best possible up the the value of the constant~$C$. 
Thus in this case the chromatic number of~$H$ is the relevant parameter which governs the
existence of perfect matchings in graphs of large minimum degree.
The first proposition deals with the case when $\chi(H)\ge 3$
as well as the case when $\chi(H)=2$ and $\hcf_\chi(H)\ge 3$.

\begin{prop}\label{bestposs1}
Let $H$ be a graph with $2\le\chi(H)=:\ell$ and let $k\in\mathbb{N}$.
Let $G_1$ be the complete $\ell$-partite graph of order $k|H|$ whose
vertex classes $U_1,\dots,U_\ell$ satisfy
$|U_1|=\lfloor k|H|/\ell\rfloor+1$, $|U_2|= \lceil k|H|/\ell\rceil-1$ and
$\lfloor k|H|/\ell\rfloor\le |U_i|\le \lceil k|H|/\ell\rceil$ for all $i\ge 3$.
(So $\delta(G_1)=\lceil(1-1/\chi(H))|G_1|\rceil-1$.)
If $\ell\ge 3$ and $\hcf_\chi(H)\neq 1$ or if $\ell=2$ and $\hcf_\chi(H)\ge 3$
then $G_1$ does not contain a perfect $H$-packing.
\end{prop}
\proof
Suppose first that $\ell\ge 3$ and that $\hcf_\chi(H)$ is finite. In this case there
are vertex classes $U_{i_1}$ and $U_{i_2}$ such that $|U_{i_1}|-|U_{i_2}|=1$
(note that these are not necessarily $U_1$ and $U_2$).
Consider any $H$-packing in $G_1$ consisting of $H_1,\dots,H_r$ say.
By induction one can show that
$|U_{i_1}\setminus(H_1\cup \dots\cup H_r)|-
|U_{i_2}\setminus(H_1\cup \dots\cup H_r)|\equiv 1\mod \hcf_\chi(H)$.
But as $\hcf_\chi(H)\neq 1$ this implies that at least one of
$U_{i_1}\setminus(H_1\cup \dots\cup H_r)$ and $U_{i_2}\setminus(H_1\cup \dots\cup H_r)$
has to be non-empty. Thus the $H$-packing $H_1,\dots,H_r$ cannot be perfect.
If $\ell\ge 3$ but $\hcf_\chi(H)=\infty$ then the colour classes of~$H$ have
the same size and thus~$G_1$ also works.

The case when $\ell=2$ is similar except that we have to work with $U_1$ and
$U_2$ and can only assume that $|U_1|-|U_2|\in\{1,2\}$.
\endproof 

In the next proposition we consider the case when $\chi(H)=2$
and $\hcf_c(H)\neq 1$. We omit its proof as it is similar to the proof
of Proposition~\ref{bestposs1}.%
     \COMMENT{Indeed, let $A$ and $B$ denote the vertex classes
of~$G_2$ such that $|A|\ge |B|$. Suppose first that $\hcf_c(H)=2$ and $k|H|$ is
not divisible by~$4$. Note that $|H|$ is even since $\hcf_c(H)=2$.
Thus $k|H|$ is divisible by $2$ (so the def of $G_2$ makes sense)
and $k|H|/2$ is an odd number. Now consider any $H$-packing in $G_2$
consisting of $H_1,\dots,H_r$ say. 
By induction one can show that
$|A\setminus(H_1\cup \dots\cup H_r)|\equiv 1\mod \hcf_c(H)$.
Thus $A\setminus(H_1\cup \dots\cup H_r)$ is non-empty and so
$H_1,\dots,H_r$ cannot be perfect.
Next consider the case when $k|H|$ is not divisible by~$2$.
Then $|A|=|B|+1$. So at least one of $|A|,|B|$ is not divisible
by~$\hcf_c(H)$. Suppose that this is true for~$A$. Then the above
argument works.
Finally, consider the case when $k|H|$ is divisible by~$2$.
Suppose first that $\hcf_c(H)\ge 3$. Then we can argue as before
since at least one of $|A|, |B|$ is not divisible by~$\hcf_c(H)$.
So suppose that $\hcf_c(H)=2$ and (as we are not in the first case)
$k|H|$ is divisible by~$4$. Then $|A|$ is odd and our argument
works again.}

\begin{prop}\label{bestposs2}
Let $H$ be a bipartite graph with $\hcf_c(H)\neq 1$ and let $k\in\mathbb{N}$.
If $\hcf_c(H)=2$ and $k|H|$ is not divisible by~$4$ let
$G_2$ be the disjoint union of two cliques of order~$k|H|/2$.
Otherwise let $G_2$ be the disjoint union of two cliques of
orders~$\lfloor k|H|/2\rfloor+1$ and~$\lceil k|H|/2\rceil-1$.
(So $\delta(G_2)\ge (1-1/\chi(H))|G_1|-2$.)
$G_2$ does not contain a perfect $H$-packing.%
\noproof
\end{prop}

The following corollary gives a characterization of those graphs with $\hcf(H)=1$.
It follows immediately from Propositions~\ref{bestposs1} and~\ref{bestposs2}
as well as Lemmas~\ref{completemove}--\ref{completebipmove}
in Section~\ref{sec:complete}. We will not need it in the proof of Theorem~\ref{thmmain}
but state it as the characterization may be of independent interest.

\begin{cor} \label{equivalence}
Let $H$ be a graph with $2\le \chi(H)=:\ell$.
Let $k'\gg |H|$ be an integer and let $G_1$ and $G_2$ be the graphs
defined in Propositions~\ref{bestposs1} and~\ref{bestposs2} for $k:=\ell k'$.
If $\chi(H)\ge 3$ then $G_1$ contains a perfect $H$-packing if and
only if $\hcf(H)=1$. Similarly, if $\chi(H)=2$ then both $G_1$ and $G_2$
contain a perfect $H$-packing if and only if $\hcf(H)=1$.%
\noproof  
\end{cor}


\section{The Regularity lemma and the Blow-up lemma}\label{sec:RLandBL}

The purpose of this section is to collect all the information we need
about the Regularity lemma and the Blow-up lemma.
See~\cite{KSi} and~\cite{JKblowup} for surveys about these.
Let us start with some more notation.
The \emph{density} of a bipartite graph $G=(A,B)$ is defined to be
$$d_G(A,B):=\frac{e_G(A,B)}{|A||B|}.$$
We also write $d(A,B)$ if this is unambiguous.
Given $\eps>0$, we say that $G$ is \emph{$\eps$-regular} if for all
sets $X\subseteq A$ and $Y\subseteq B$ with $|X|\ge \eps |A|$ and
$|Y|\ge \eps |B|$ we have $|d(A,B)-d(X,Y)|<\eps$. Given $d\in[0,1]$, we say
that $G$ is \emph{$(\eps,d)$-superregular} if all
sets $X\subseteq A$ and $Y\subseteq B$ with $|X|\ge \eps |A|$ and
$|Y|\ge \eps |B|$ satisfy $d(X,Y)>d$ and, furthermore, if $d_G(a)>d|B|$
for all $a\in A$ and $d_G(b)> d|A|$ for all $b\in B$.

We will use the
following degree form of Szemer\'edi's Regularity lemma which can be easily
derived from the classical version. Proofs of
the latter are for example included in~\cite{BGraphTh} and~\cite{Diestel}.

\begin{lemma}[Regularity lemma]\label{deg-reglemma}
For all $\eps>0$ and all integers $k_0$ there is an $N=N(\eps,k_0)$ such
that for every number $d\in [0,1]$ and for every graph $G$ on at least
$N$ vertices there exist a
partition of $V(G)$ into $V_0,V_1,\dots,V_k$ and a spanning subgraph
$G'$ of $G$ such that the following holds:
\begin{itemize}
\item $k_0\le k\le N$,
\item $|V_0|\le \eps |G|$,
\item $|V_1|=\dots=|V_k|=:L$,
\item $d_{G'}(x)>d_G(x)-(d+\eps)|G|$ for all vertices $x\in G$,
\item for all $i\ge 1$ the graph $G'[V_i]$ is empty,
\item for all $1\le i<j\le k$ the graph $(V_i,V_j)_{G'}$ is $\eps$-regular
and has density either $0$ or $>d$.
\end{itemize}
\end{lemma}

The sets $V_i$ ($i\ge 1$) are called \emph{clusters}, $V_0$ is called
the \emph{exceptional set}.
Given clusters and $G'$ as in Lemma~\ref{deg-reglemma}, the
\emph{reduced graph} $R$ is the graph whose vertices are
$V_1,\dots,V_k$ and in which $V_i$ is joined to
$V_j$ whenever $(V_i,V_j)_{G'}$ is $\eps$-regular and has density $>d$. Thus
$V_iV_j$ is an edge of $R$ if and only if $G'$ has an edge between $V_i$
and $V_j$.

Given a set $A\subseteq V(R)$, we call the set of all those vertices of
$G$ which are contained in clusters belonging to $A$ the \emph{blow-up of $A$}.
Similarly, if $R'$ is a subgraph of $R$, then the \emph{blow-up of $R'$}
is the subgraph of $G'$ induced by the blow-up of~$V(R')$.

We will also use the Blow-up lemma of Koml\'os, S\'ark\"ozy and
Szemer\'edi~\cite{KSSblowup}. It implies that dense regular pairs behave
like complete bipartite graphs with respect to containing bounded degree
graphs as subgraphs.

\begin{lemma}[Blow-up lemma]\label{blowup}
Given a graph $F$ on $\{1,\dots,f\}$ and positive numbers $d,\Delta$, there is
a positive number $\eps_0=\eps_0(d,\Delta,f)$
such that the following holds. Given
$L_1,\dots,L_f\in \mathbb{N}$ and $\eps\le \eps_0$, let $F^*$ be the graph obtained
from $F$ by replacing each vertex $i\in F$ with a set $V_i$ of $L_i$ new
vertices and joining
all vertices in $V_i$ to all vertices in $V_j$ whenever $ij$ is an edge
of $F$. Let $G$ be a spanning subgraph of $F^*$
such that for every edge $ij\in F$ the graph $(V_i,V_j)_G$ is
$(\eps,d)$-superregular. Then $G$ contains a copy of every subgraph $H$ of
$F^*$ with $\Delta(H)\le \Delta$.
\end{lemma}

\section{Preliminaries and overview of the proof}\label{sec:overview}

Let $H$ be a graph of chromatic number $\ell\ge 2$.
Put
\begin{equation}\label{eqdefxi}
z_1:=(\ell-1)\sigma(H),\ \  z:=|H|-\sigma(H),\ \ 
\xi:=\frac{z_1}{z}=\frac{(\ell-1)\sigma(H)}{|H|-\sigma(H)}.
\end{equation}
(Recall that $\sigma(H)$ is the smallest colour class in any $\ell$-colouring
of $H$.) Note that $\xi<1$ if $\chi_{cr}(H)<\chi(H)$ and so in particular if $\hcf(H)=1$.
Let $B^*$ denote the complete $\ell$-partite graph with one vertex class of size
$z_1$ and $\ell-1$ vertex classes of size $z$. Note that $B^*$ has a perfect
$H$-packing consisting of $\ell-1$ copies of $H$. Moreover,
it is easy to check that
\begin{equation}\label{eqchicr}
\chi_{cr}(H)=\chi_{cr}(B^*)=\ell-1+\xi.
\end{equation}
Call $B^*$ the \emph{bottlegraph assigned to~$H$}. We now give an overview of the
proof of Theorem~\ref{thmmain}. In Section~\ref{sec:applyRG} we first apply
the Regularity lemma to $G$ in order to obtain a set $V_0$ of exceptional vertices
and a reduced graph~$R$. It will turn out that the minimum degree of $R$ is
almost $(1-1/\chi_{cr}(B^*))|R|$. So $R$ has an almost perfect $B^*$-packing $\B'$
by Theorem~\ref{thmKomlos}. Let $B_1,\dots,B_{k'}$
denote the copies of $B^*$ in~$\B'$. Our aim in Sections~\ref{sec:adjust}
and~\ref{sec:div} is to show that one can take out a small number of
suitably chosen copies of $B^*$ from $G$ to achieve that the following conditions
hold:
\begin{itemize}
\item[($\alpha$)] Each vertex in $V_0$ lies in one of these copies of $B^*$ taken out
from~$G$. Moreover, each vertex that does not belong to a blow-up of some
$B_t\in \B'$ also lies in one of these copies of~$B^*$. 
\item[($\beta$)] The (modified) blow-up of each $B_t\in \B'$ has a 
perfect $H$-packing.
\end{itemize}
Note that if we say that we take out a copy of $B^*$ (or $H$) from $G$ then we mean that
we delete all its vertices from $G$ and thus also from the clusters they belong to.
So in particular the (modified) blow-up of $B_t$ no longer contains these vertices.

For all $t\le k'$ and all $j\le \ell$ let $X_j(t)$ denote the (modified) blow-up of the $j$th vertex
class of $B_t\in \B'$ (where the $\ell$th vertex class of $B_t$ is the small one).
It will turn out that~$(\beta)$ holds if the $X_j(t)$ satisfy the conditions
in the following definition. (The graph $G'\subseteq G$ obtained from the
Regularity lemma will play the role of $G^*$ in Definition~\ref{defblownupcover}.
So it will be easy to satisfy condition~(a) of Definition~\ref{defblownupcover}.)

\begin{defin}\label{defblownupcover}
{\rm Suppose that $G$ is a graph whose order~$n$ is divisible by~$|B^*|$.
Let $k\ge 1$ be an integer and let $\eps\ll d\ll\beta\ll 1$ be positive
constants.
We say that \emph{$G$ has a blown-up $B^*$-cover for parameters
$\eps,d,\beta,k$} if there exists a spanning subgraph $G^*$ of $G$
and a partition $X_1(1),\dots,X_\ell(1),\dots,X_1(k),\dots,X_\ell(k)$
of the vertex set of $G$ such that the following holds:
\begin{itemize}
\item[(a)] All the bipartite subgraphs $(X_j(t),X_{j'}(t))_{G^*}$ of $G^*$
between $X_j(t)$ and $X_{j'}(t)$ are $(\eps,d)$-superregular
whenever $j\neq j'$.
\item[(b)] $|X_1(t)\cup\dots\cup X_\ell(t)|$
is divisible by~$|B^*|$ for all $t\le k$.
\item[(c)] $(1-\beta^{1/10})|X_\ell(t)|\le \xi |X_j(t)|\le (1-\beta)|X_\ell(t)|$
for all $j<\ell$ and all
$t\le k$ and $|\, |X_j(t)|-|X_{j'}(t)|\,|\le d|X_1(t)\cup\dots\cup X_\ell(t)|$ for all $1\le j<j'<\ell$
and all $t\le k$.
\end{itemize}}
\end{defin}

For all $t \le k$, we call the $\ell$-partite subgraph of $G^*$ whose vertex
classes are the sets $X_1(t),\dots,X_\ell(t)$
the \emph{$t$'th element} of the blown-up $B^*$-cover. 
The \emph{complete $\ell$-partite graph corresponding to the $t$'th element}
is the one whose vertex classes have sizes $|X_1(t)|,\dots,|X_\ell(t)|$.
Note that condition~(c) implies that for all $j<\ell$ the ratio of $|X_j(t)|$
to the total size of the $t$th element is a little smaller than $z/|B^*|$
(recall that $z$ is the size of the large vertex classes of~$B^*$).

The following lemma implies that the complete $\ell$-partite graph corresponding to
some element of a blown-up $B^*$-cover contains
a perfect $H$-packing. Combined with the Blow-up lemma, this will imply that each 
element of the blown-up $B^*$-cover has a perfect $H$-packing.
(Thus~$(\beta)$ will be satisfied if the $X_j(t)$ are as in
Definition~\ref{defblownupcover}.)

\begin{lemma}\label{corcomplete1}
Let $H$ be a graph with $\ell:=\chi(H)\ge 2$ and
$\hcf(H)=1$. Let $\xi$ be as defined in~$(\ref{eqdefxi})$.
Let $0<d\ll\beta\ll \xi,1-\xi,1/|H|$ be positive constants. 
Suppose that $F$ is a complete $\ell$-partite graph with vertex classes
$U_1,\dots,U_\ell$ such that $|F|\gg |H|$ is divisible by $|H|$,
$(1-\beta^{1/10})|U_\ell|\le \xi |U_i|\le (1-\beta)|U_\ell|$ for all $i<\ell$ and such that
$|\, |U_i|-|U_j|\,|\le d|F|$ whenever $1\le i<j<\ell$. Then $F$ contains a
perfect $H$-packing.
\end{lemma}

Lemma~\ref{corcomplete1} will be proved at the end of Section~\ref{sec:complete}, where we will
deduce it from Lemmas~\ref{completeconst}
and~\ref{completeapprox}, which are also proved in that section.
Lemma~\ref{corcomplete1} is one of the points where the condition that $\hcf(H)=1$ is necessary.

The following lemma shows that we can find a blown-up $B^*$-cover as long as $G$
satisfies certain properties. Roughly speaking these properties~(i) and~(ii) say
that $G$ is not too close to being one of the extremal graphs having minimum degree
almost $(1-1/\chi_{cr}(H))n$ but not containing a perfect $H$-packing.
We will refer to this as the non-extremal case.
 
\begin{lemma}\label{nonextremal}
Let $H$ be a graph of chromatic number $\ell\ge 2$ such that $\chi_{cr}(H)<\ell$.
Let $B^*$ denote the bottlegraph assigned to~$H$ and let $z$ and $\xi$
be as defined in~$(\ref{eqdefxi})$. Let
$$
\eps'\ll d'\ll\theta\ll \tau\ll\xi,1-\xi,1/|B^*|
$$
be positive constants.
There exist integers $n_0$ and $k_1=k_1(\eps',\theta,B^*)$ such that
the following holds.
Suppose that $G$ is a graph whose order $n\ge n_0$ is divisible by
$|B^*|$ and whose minimum degree satisfies
$\delta(G)\ge (1-\frac{1}{\chi_{cr}(H)}-\theta)n$. Furthermore,
suppose that $G$ satisfies the following further properties:
\begin{itemize}
\item[{\rm (i)}] $G$ does not contain a
vertex set $A$ of size $z n/|B^*|$ such that $d(A)\le \tau$.
\item[{\rm (ii)}] Additionally, if $\ell=2$ then $G$ does not
contain a vertex set $A$ such that $d(A,V(G)\setminus A)\le \tau$.
\end{itemize}
Then there exists a family $\B^*$ of at most $\theta^{1/3}n$ disjoint copies
of $B^*$ in $G$ such that the graph $G-\bigcup \B^*$ (which is obtained from $G$
by taking out all the copies of $B^*$ in $\B^*$) has a blown-up $B^*$-cover
with parameters~$2\eps',d'/2,2\theta,k_1$.
\end{lemma}
Lemma~\ref{nonextremal} will be proved in Section~\ref{sec:nonextremal}.
Lemmas~\ref{corcomplete1} and~\ref{nonextremal} together with the Blow-up lemma
imply that in the non-extremal case we can satisfy conditions~($\alpha$)
and~($\beta$), i.e.~we have a perfect $H$-packing in this case. 
This is formalized in the following corollary.

\begin{cor}\label{cornonextremal}
Let $H$ be a graph of chromatic number $\ell\ge 2$ such that
$\hcf(H)=1$. Let $B^*$ denote the bottlegraph assigned
to~$H$ and let $z$ and $\xi$ be as defined in~$(\ref{eqdefxi})$. Let
$\theta\ll\tau\ll\xi,1-\xi,1/|B^*|$ be positive constants.
There exists an integer $n_0$ such that the following holds.
Suppose that $G$ is a graph whose order $n\ge n_0$ is  divisible by~$|B^*|$
and whose minimum degree satisfies
$\delta(G)\ge (1-\frac{1}{\chi_{cr}(H)}-\theta)n$. Furthermore,
suppose that $G$ satisfies the following further properties:
\begin{itemize}
\item[{\rm (i)}] $G$ does not contain a
vertex set $A$ of size $z n/|B^*|$ such that $d(A)\le \tau$.
\item[{\rm (ii)}] Additionally, if $\ell=2$ then $G$ does not
contain a vertex set $A$ such that $d(A,G-A)\le \tau$.
\end{itemize}
Then $G$ has a perfect $H$-packing.
\end{cor}

\removelastskip\penalty55\medskip\noindent{\bf Proof of
Corollary~\ref{cornonextremal}. }
Fix positive constants $\eps'$, $d'$ such that
$$
\eps'\ll d'\ll\theta\ll \tau\ll \xi,1-\xi,1/|B^*|.
$$
An application of Lemma~\ref{nonextremal} shows that by taking out a small number
of disjoint copies of $B^*$ from $G$ we obtain a subgraph which has a blown-up
$B^*$-cover with parameters $2\eps',d'/2,2\theta,k_1$.
Conditions~(b) and (c) in Definition~\ref{defblownupcover}
imply that the complete $\ell$-partite graphs corresponding to the $k_1$ elements of
this blown-up $B^*$-cover satisfy the assumptions of Lemma~\ref{corcomplete1}
with $d:=d'/2$ and where $2\theta$ plays the role of~$\beta$.
Thus each of these complete $\ell$-partite graphs contains a perfect $H$-packing.
Condition~(a) in Definition~\ref{defblownupcover} ensures that we can now apply the 
Blow-up lemma (Lemma~\ref{blowup}) to each of the $k_1$ elements in the blown-up
$B^*$-cover to obtain a perfect $H$-packing of this element. All these $H$-packings
together with the copies of $B^*$ taken out earlier in order to obtain the blown-up
$B^*$-cover yield a perfect $H$-packing of~$G$.
\endproof

The extremal cases (i.e.~where $G$ satisfies either~(i) or~(ii)) will be dealt with in
Section~\ref{sec:extremal}. These cases also rely on Lemma~\ref{nonextremal}.
For example if $G$ satisfies (i)
but $G-A$ does not satisfy (i) or (ii), then very roughly the strategy is to apply
Lemma~\ref{nonextremal} to find a perfect $B_1^*$-packing of $G-A$,
where $B_1^*$ is obtained from $B^*$ by removing one of the large colour classes.
The minimum degree of $G$ will ensure that the bipartite subgraph
spanned by $A$ and $V(G)-A$ is almost complete. This will be used to extend the $B_1^*$-packing
of $G-A$ to a perfect $B^*$-packing of~$G$. The reason we considered sets $A$ of size
$zn/|B^*|$ in~(i) is that this is precisely the number of vertices needed to extend each copy
of $B^*_1$ to a copy of~$B^*$. (Recall that $z$ was the size of the large vertex classes
of~$B^*$.)

However, for this strategy to work, we first need to modify the set $A$ slighly.
We will also take out some carefully chosen copies of $H$ from $G$.
One matter which complicates the argument is that
$B_1^*$ does not necessarily satisfy $\hcf(B_1^*)=1$. 
This means that we cannot find perfect $B_1^*$-packing of $G-A$
by a direct application of Lemma~\ref{nonextremal}, as the blown-up $B_1^*$-cover produced by that lemma
does not necessarily yield a perfect $B_1^*$-packing of $G-A$. To overcome this difficulty,
we will work directly with the blown-up $B_1^*$-cover. 
So the use of Lemma~\ref{nonextremal}
in Section~\ref{sec:extremal} is the reason why we do not assume 
$\hcf(H)=1$ in Lemma~\ref{nonextremal}. It is also the reason why we
allow for an error term $\theta n$ in the minimum degree condition on $G$.

\section{The non-extremal case: proof of Lemma~\ref{nonextremal}}\label{sec:nonextremal}

The purpose of this section is to prove Lemma~\ref{nonextremal}.

\subsection{Applying the Regularity lemma and choosing a packing of the
reduced graph}\label{sec:applyRG}
We will fix further constants satisfying the following hierarchy
\begin{equation}\label{eqconst}
0< \eps\ll \eps'\ll d'\ll d\ll \theta\ll\tau\ll\xi, 1-\xi,1/|B^*|.
\end{equation}
Moreover, we choose an integer $k_0$ such that
\begin{equation}\label{eqk0}
k_0\ge n_1(\theta,B^*),
\end{equation}
where $n_1$ is as defined in Theorem~\ref{thmKomlos}.
We put
\begin{equation}\label{eqdefk1}
k_1:=\lfloor N(\eps,k_0)/|B^*|\rfloor,
\end{equation}
where $N(\eps,k_0)$ is as defined in the Regularity lemma (Lemma~\ref{deg-reglemma}).
In what follows, we assume that the order $n$ of our given graph $G$ is sufficiently
large for our estimates to hold.
We now apply the Regularity lemma with parameters $\eps$, $d$ and $k_0$ to $G$
to obtain clusters, an exceptional set $V_0$, a spanning subgraph $G'\subseteq G$
and a reduced graph $R$. (\ref{eqconst}) together with the well-known fact that the
minimum degree of~$G$ is almost inherited by its reduced graph (see
e.g.~\cite[Prop.~9]{KOTplanar} for an explicit proof) implies that
\begin{align}\label{eqminR1}
\delta(R) \ge \left(1-\frac{1}{\chi_{cr}(H)}-2\theta\right)|R|
\stackrel{(\ref{eqchicr})}{=}\left(1-\frac{1}{\ell-1+\xi}-2\theta\right)|R|.
\end{align}
Since $|R|\ge k_0\ge n_1(\theta,B^*)$ by~(\ref{eqk0}),
we may apply Theorem~\ref{thmKomlos} to $R$ to find a
$B^*$-packing $\B'$ which covers all but at most $\sqrt{\theta}|R|$ vertices of~$R$.
(More precisely, we apply Theorem~\ref{thmKomlos} to a graph $R'$ which is obtained from $R$
by adding at most $\theta^{3/4}|R|$ new vertices and connecting them to all other vertices.
By~(\ref{eqchicr}) and~(\ref{eqminR1}), we have
$\delta(R')\ge (1-1/\chi_{cr}(B^*))|R'|$, as required in Theorem~\ref{thmKomlos}.%
     \COMMENT{Note that it does not suffice to add $2\theta |R|$ vs. Indeed, let $y$ denote
the number of vertices which we have to add. Then $y$ has to satisfy
$$(1-\frac{1}{\ell-1+\xi}-2\theta)|R|+y\ge (1-\frac{1}{\ell-1+\xi})(|R|+y)$$
and thus $y\ge 2\theta|R|(\ell-1+\xi)$.}
Removing the new vertices results in a $B^*$-packing of $R$ which has the desired size.) 
We delete all the clusters not contained in some copy of $B^*$ in $\B'$
from $R$ and add all the vertices lying in these clusters to the exceptional set~$V_0$.
Thus $|V_0|\le \eps n+\sqrt{\theta} n\le 2\sqrt{\theta} n$.
From now on, we denote by $R$ the subgraph
of the reduced graph induced by all the remaining clusters. Thus $\B'$ now is a perfect
$B^*$-packing of~$R$ and we still have that
\begin{align}\label{eqmindegreduced}
\delta(R) \ge \left(1-\frac{1}{\ell-1+\xi}-2\sqrt{\theta}\right)|R|.
\end{align}
It is easy to check that for all $B\in\B'$ we can
replace each cluster $V_a$ in $B$ by a subcluster of size $L':=(1-\eps|B^*|)L$ such that for
each edge $V_aV_b$ of $B$ the bipartite subgraph of $G'$ between the chosen
subclusters of $V_a$ and $V_b$ is $(2\eps,d/2)$-superregular
(see e.g.~\cite[Prop.~8]{KOTplanar}).
Add all the vertices of~$G$ which do not lie in one
of the chosen subclusters to the exceptional set $V_0$. Then
\begin{equation*}
|V_0|\le 3\sqrt{\theta} n.
\end{equation*}
By adjusting $L'$ if necessary and adding a bounded number of
further vertices to~$V_0$ we may assume that $L'$ is divisible by $z_1z$.
(Recall that $z_1$ and $z$ were defined in~(\ref{eqdefxi}).)
From now on, we refer to the chosen subclusters as the clusters of~$R$.

Next we partition
each of these clusters $V_a$ into a red part $V_a^{red}$ and a blue part $V_a^{blue}$
such that $|V_a^{red}|=\theta^2 |V_a|$ and such that
$|\, |N_G(x)\cap V_a^{red}|-\theta^2|N_G(x)\cap V_a|\, |\le \eps L'$ for every vertex $x\in G$.
(Consider a random partition to see that there are $V_a^{red}$ and $V_a^{blue}$
with these properties.) Together all these partitions of the clusters of $R$ yield a
partition of the vertices of $G-V_0$ into red and blue vertices.
We will use these partitions to ensure that even after some modifications which we have
to carry out during the proof, the edges
of the $B\in\B'$ will still correspond to superregular subgraphs of~$G'$.
More precisely, during the proof we will take out certain copies of $B^*$ from $G$,
but each copy will avoid all the red vertices. All the vertices contained in
these copies of $B^*$ will be removed from the clusters they belong to.
However, if we look at the (modified) bipartite subgraph of $G'$ which corresponds to some edge
$V_aV_b$ of $B\in \B'$, then this subgraph of $G'$ will still be $(\eps', d')$-superregular
since it still contains all vertices in $V_a^{red}$ and $V_b^{red}$.

The blown-up $B^*$-cover required in Lemma~\ref{nonextremal} will
be otained from the cover corresponding to $\B'$ by taking out a small number of
copies of~$B^*$ from~$G$. This will be done in two steps. Firstly, we will
take out copies of $B^*$ to ensure that for every $B\in \B'$ the size of
the blow-up of each of its $\ell-1$ large vertex classes is significantly smaller than
$(z/|B^*|)$times the size of the blow-up of the entire~$B$.
This will ensure that the cover corresponding to $\B'$ satisfies condition~(c) in the definition of
blown-up $B^*$-cover (Definition~\ref{defblownupcover}).
Moreover, each exceptional vertex will be contained in one of the copies
taken out. So after this step we also have incorporated all the exceptional
vertices. All the copies of $B^*$
deleted in this process will avoid the red vertices of~$G$.
In the second step we will then take out a bounded number of further copies of
$B^*$ in order to achieve that the blow-up of each $B\in\B'$ is divisible by~$|B^*|$
(as required in condition~(b) in Definition~\ref{defblownupcover}).
As mentioned in the previous
paragraph, the blown-up $B^*$-cover thus obtained from $\B'$ will also satisfy
condition~(a) in Definition~\ref{defblownupcover} since in this last step we
remove only a bounded number of further vertices from the clusters, which does not affect
the superregularity significantly.
Finally, since the blown-up $B^*$ cover obtained in this way has $|\B'|\le k_1$
elements, we may have to split some of the blown up copies of $B^*$
to obtain a blown-up $B^*$ cover with exactly $k_1$ elements.

\subsection{Adjusting the sizes of the vertex classes in the blow-ups
of the $B\in \B'$}\label{sec:adjust}

Let $B_1,\dots,B_{k'}$ denote the copies of $B^*$ in~$\B'$. 
As described at the end of the previous section,
our next aim is to take out a small number of copies of $B^*$ from $G$ to achieve that,
for all $t\le k'$, the blow-up of each large vertex class of $B_t$ is
significantly smaller than $(z/|B^*|)$times the size of the blow-up of~$B_t$ itself.
It turns out that this becomes simpler if we first split the blow-up of each $B_t$
into $z_1z$ `smaller blow-ups'. Then we take out copies of $B^*$ from $G$
in order to modify the sizes of these smaller blow-ups. This will imply that
sizes of the blown-up vertex classes of the original $B_t$'s are as desired.
We will not remove red vertices in this process.

Thus consider any $B_t$. We will think of the $\ell$th vertex class of~$B_t$
as the one having size~$z_1$. For all $j<\ell$,
split each of the $z$ clusters belonging to the $j$th vertex class
of $B_t$ into $z_1$ subclusters of equal size.
Let $Z'_j((t-1)z_1z+1),\dots, Z'_j(tz_1z)$ denote the subclusters thus obtained.
Similarly, split each cluster belonging to the $\ell$th vertex class of $B_t$ into
$z$ subclusters of equal size. Let $Z_\ell((t-1)z_1z+1),\dots, Z_\ell(tz_1z)$
denote the subclusters thus obtained. 
Put
$$k'':=z_1zk'.
$$
Given $i\le k''$, we think of the $\ell$-partite subgraph of
$G'$ with vertex classes $Z'_1(i),\dots,Z'_{\ell-1}(i),Z_\ell(i)$ as a
blown-up copy of $B^*$. (Indeed, note that $\xi |Z'_j(i)|=|Z_\ell(i)|$
for all $j<\ell$.) We may assume that about $\theta^2|Z'_j(i)|$ vertices in $|Z'_j(i)|$
are red and that for every vertex $x\in G$ about a $\theta^2$-fraction of its
neighbours in each~$Z'_j(i)$ are red
and that the analogue holds for~$Z_\ell(i)$.
(Indeed, consider random partitions
again to show that this can be guaranteed.)

In order to achieve that the size of each large vertex class
is significantly smaller than the size of the entire blown-up copy,
we will remove a $\theta$-fraction of vertices from each of
$Z'_1(i),\dots,Z'_{\ell-1}(i)$.
We will add all these vertices to~$V_0$. The aim then is to incorporate all
the vertices in~$V_0$ by taking out copies of $B^*$ from~$G$.
Of course, this has to be done in such a way that we don't destroy the
properties of the vertex classes again. So put
$$L'':=(1-\theta)L'/z_1.
$$
For all $i\le k''$ and all $j< \ell$ remove $\theta L'/z_1$ blue vertices
from $Z'_j(i)$ and add all these vertices to~$V_0$.
Denote the subset of $Z'_j(i)$ thus obtained by~$Z_j(i)$.
So
\begin{equation}\label{sizeZji}
|Z_\ell(i)|=L'/z=\xi L''/(1-\theta),\ \ \ |Z_j(i)|=L''=(1-\theta)|Z_\ell(i)|/\xi
\end{equation}
whenever $j<\ell$. Also, we now have that
\begin{equation}\label{eqV0}
|V_0|\le 4\sqrt{\theta}n.
\end{equation}
Note that for all $i\le k''$ and all
$0\le j<j'\le \ell$ the graph $(Z_j(i),Z_{j'}(i))_{G'}$ is $\eps'$-regular
and has density at least~$d'$.
Denote by $Z^{red}_j(i)$ the set of red vertices in~$Z_j(i)$.

Let us now prove the following claim. Roughly speaking, it states that by taking out
a small number of copies of $B^*$ from $G$ we can incorporate all the
vertices in~$V_0$ and that this can be done without destroying the properties
of the vertex classes of the blown-up copies of~$B^*$.

\medskip

\noindent
\textbf{Claim.} {\it We can take out at most $3|V_0|\le 12\sqrt{\theta}n$
disjoint copies of $B^*$ from $G$ which cover all the vertices in~$V_0$ and have the property
that the leftover sets $Y_j(i)\subseteq Z_j(i)$ thus obtained
satisfy
\begin{itemize}
\item[{\rm (a)}] $Z^{red}_j(i)\subseteq Y_j(i)$,
\item[{\rm (b)}] $L''-|Y_1(i)|\le \theta^{1/7}L''$,
\item[{\rm (c)}] $|Y_1(i)|=\dots=|Y_{\ell-1}(i)|\le (1-\theta)|Y_\ell(i)|/\xi$.
\end{itemize}}

\medskip

\noindent
To prove this claim, we show that for every vertex $x\in V_0$ in turn
we can take out either one, two or three disjoint copies of $B^*$
which satisfy the following three properties.%
   \textno Firstly, $x$ lies in one of
the copies. Secondly, these copies avoid all the red vertices.
Thirdly, when removing these copies from $G$ then,
for every $i\le k''$, we either delete no vertex at all in
$Z_1(i)\cup\dots\cup Z_\ell(i)$ or else we delete precisely $z$ vertices
in each of $Z_1(i),\dots,Z_{\ell-1}(i)$ and delete either $z_1$ or $z_1-1$
vertices in~$Z_\ell(i)$. &(*)%

\noindent
Together with~(\ref{sizeZji}) this implies that
after each step the subsets obtained
from the~$Z_j(i)$ will satisfy conditions~(a) and (c).%
     \COMMENT{Indeed, for (c) we need that
$\xi(|Z_1(i)|-z)\le (1-\theta)(|Z_\ell(i)|-z_1)$. But (\ref{sizeZji}) together
with the fact that $-z_1=-\xi z\le -(1-\theta)z_1$ show that this holds.}
We will discuss later how~(b) can be satisfied too. 

Thus consider the first vertex~$x\in V_0$. To find the copies of $B^*$
satisfying~$(*)$ we will distinguish several cases.
Suppose first that there exists
an index $i=i(x)$ such that $x$ has at least $\theta L''$ neighbours
in $Z_j(i)$ for all $j<\ell$.%
     \footnote{In later steps we will ask whether a vertex $x'\in V_0$ has
at least $\theta L''$ neighbours in the \emph{current} set $Z_j(i)$
for all $j<\ell$.}
Take out a copy of $B^*$ from $G$
which contains~$x$, which meets each of $Z_1(i),\dots,Z_{\ell-1}(i)$ in $z$
vertices and $Z_\ell(i)$ in $z_1-1$ vertices and which avoids the red
vertices of~$G$. (The existence of such copies of
$B^*$ in $G$ easily follows from a `greedy' argument based on the 
$\eps'$-regularity of the bipartite subgraphs $(Z_j(i),Z_{j'}(i))_{G'}$ of $G'$,
see e.g.~Lemma~7.5.2 in~\cite{Diestel}
or Theorem~2.1 in~\cite{KSi}. We will often use this and similar facts below.
We can avoid the red vertices since $|Z_j^{red}(i)|\ll \theta L''$ and so
most of the neighbours of $x$ in $Z_j(i)$ will be blue.)%
     \COMMENT{Indeed, $|Z_j^{red}(i)|\approx \theta^2L'/z_1=\theta^2(1-\theta)L''
\ll \theta L''$.}

Next suppose that we cannot find an index $i$ as above.
By relabelling if necessary, we may assume that $x$ has
at most $\theta L''$ neighbours in $Z_1(i)$ for all $i\le k''$.
Let $I$ denote the set of all those indices $i$ for which
$x$ has at least $\theta L''$ neighbours in $Z_j(i)$
for all $j=2,\dots,\ell$. To obtain a lower bound on the
size of~$I$, we now consider~$d_{G'}(x)$. This shows that 
\begin{align*}
(k''-|I|) L''(\ell-2+2\theta)+|I|L''(\ell-2+\xi/(1-\theta)+\theta)
& \stackrel{(\ref{sizeZji})}{\ge} d_{G'}(x)-|V_0|\\
& \stackrel{(\ref{eqV0})}{\ge} (1-1/\chi_{cr}(H)-\theta^{1/3})n.
\end{align*}
The term $\theta^{1/3}n$ in the second line is a bound on $|V_0|$ (with room to
spare -- this  will be useful later on).
The above equation implies that%
     \COMMENT{Note that the $\theta^{1/3}n$ also works
in later steps since we claim that we remove at most $12\sqrt{\theta}n$
copies of $B^*$ in the entire process and since there are $\ll \theta^{1/3}k''$
indices $i$ which we exclude from consideration because the size of some
$Z_j(i)$ is critical. Thus the modified inequality looks
like
\begin{align*}
\text{old LHS} \ge & \text{degree of } x\text{ in current subgraph of } G'\\
 & -|V_0|
 -\text{ no. vs already removed }
 -\bigcup_{\text{excluded }i}\bigcup_{j=1}^\ell Z_j(i)
\ge \text{old RHS}.
\end{align*}
Moreover, we have $\xi/(1-\theta)$ instead of $\xi$ since the
size of the $Z_\ell(i)$ is $\xi L''/(1-\theta)$.
From the above inequality,
we get
\begin{align*}
\frac{\xi|I| L''}{1-\theta} & \ge \left(1-\frac{1}{\ell-1+\xi}-2\theta^{1/3}\right)n
-k''L''(\ell-2)\\
& \ge \frac{\ell-2+\xi}{\ell-1+\xi}(\ell-1+\xi)L''k''
-k''L''(\ell-2)-2\theta^{1/3}\ell k''L''\\
& = \xi L''k''-2\theta^{1/3}\ell k''L''.
\end{align*}
This shows that
$|I|\ge (1-\theta)k''-2\theta^{1/3}\ell k''\ge (1-\theta^{1/4})k''$.}
\begin{equation}\label{eqsizeI}
|I|\ge (1-\theta^{1/4})k''.
\end{equation}
Now suppose  that there are two indices $i_1,i_2\in I$ such that the density
$d_{G'}(Z_{1}(i_1),Z_j(i_2))$ is nonzero for all $j=1,\dots,\ell-1$. (Note that
the last condition in Lemma~\ref{deg-reglemma} implies that then
each of the bipartite subgraphs $(Z_{1}(i_1),Z_j(i_2))_{G'}$
of $G'$ is $\eps'$-regular and has density at least~$d'$.) In this case
we take out two disjoint copies of $B^*$. The first contains~$x$ and
has $z$ vertices in each of $Z_2(i_1),\dots,Z_{\ell-1}(i_1)$,
$z-1$ vertices in $Z_{1}(i_1)$ and $z_1$ vertices in $Z_\ell(i_1)$.
Such a copy of $B^*$ exists since $i_1\in I$.
The second copy will have one vertex in $Z_1(i_1)$, $z$ vertices in
each of $Z_1(i_2),\dots,Z_{\ell-1}(i_2)$ and $z_1-1$ vertices in
$Z_\ell(i_2)$.  Again, these copies of $B^*$ are chosen such that they
avoid the red vertices of~$G$.
So we may assume that there are no indices $i_1,i_2$ as above.

Suppose next that there are indices $i_3,i_4\in I$
and $j^*=j^*(i_3,i_4)$ with $2\le j^*\le \ell-1$ and such that
$d_{G'}(Z_{j^*}(i_3),Z_j(i_4))>0$ for all $j=2,\dots,\ell$
and $d_{G'}(Z_{j}(i_3),Z_1(i_4))>0$ for all $j^*\neq j\le \ell$.
In this case we take out 2 disjoint copies of $B^*$ again.
The first one contains~$x$, has $z$ vertices in $Z_{j^*}(i_3)$
and in each of $Z_2(i_4),\dots,Z_{\ell-1}(i_4)$ and $z_1-1$ vertices
in $Z_\ell(i_4)$. The second copy will have $z_1$ vertices in $Z_\ell(i_3)$
and $z$ vertices in $Z_1(i_4)$
as well as $z$ vertices in each $Z_j(i_3)$ with $j^*\neq j<\ell$.
Again, all these copies of $B^*$ are chosen such
that they avoid all the red vertices.
So we may assume that there are no such indices~$i_3,i_4,j^*$.

Suppose next that there are indices $i_5,i_6, i_7\in I$
and $j^\diamond=j^\diamond(i_5,i_6,i_7)$ with $2\le j^\diamond\le \ell-1$ and such that
$d_{G'}(Z_{j^\diamond}(i_5),Z_j(i_7))>0$ for all $j=1,\dots,\ell-1$
and $d_{G'}(Z_{j}(i_5),Z_1(i_6))>0$ for all $j^\diamond\neq j\le \ell$.
In this case we take out 3 disjoint copies of $B^*$.
The first copy contains~$x$, has $z-1$ vertices in $Z_{1}(i_6)$,
$z$ vertices in each of $Z_2(i_6),\dots,Z_{\ell-1}(i_6)$ and
$z_1$ vertices in~$Z_\ell(i_6)$. The second copy has one vertex
in $Z_1(i_6)$, $z-1$ vertices in $Z_{j^\diamond}(i_5)$,
$z$ vertices in each $Z_j(i_5)$ with $j^\diamond\neq j<\ell$ and $z_1$ vertices
in~$Z_\ell(i_5)$. The third copy has one vertex in $Z_{j^\diamond}(i_5)$,
$z$ vertices in each of $Z_1(i_7),\dots,Z_{\ell-1}(i_7)$ and $z_1-1$ vertices
in~$Z_\ell(i_7)$. Again, all these copies of $B^*$ are chosen such
that they avoid all the red vertices.
So we may assume that there are no such indices~$i_5,i_6,i_7,j^\diamond$.

We will show that together with our previous three assumptions this leads to
a contradiction to our assumption on the minimum degree of~$G$.
For this, first note that there are at least
$\tau |I|^2/4$ ordered pairs of indices $i,i'\in I$
for which $d_{G'}(Z_1(i),Z_1(i'))>0$. Indeed, otherwise the union $U$ of all
the $Z_1(i)$ with $i\in I$ would have density at most $\tau/2$ in $G'$ and thus
density at most $3\tau/4$ in $G$. But $zn/|B^*|-\theta^{1/10} n\le |U|\le zn/|B^*|$
by~(\ref{eqsizeI}).%
     \COMMENT{Condition~(b) implies that this also holds in later steps.}
Thus by adding at most $\theta^{1/10}n\ll \tau |U|$ vertices to $U$ if
necessary we would obtain a set
$A$ that contradicts condition~(i) of Lemma~\ref{nonextremal}.

Given $i'\in I$, we call an index $i\in I$ \emph{useful for $i'$}
if $d_{G'}(Z_1(i),Z_1(i'))>0$. 
Let $I'\subseteq I$ be the set of all those indices $i'\in I$ for
which at least  $\tau |I|/8$ other indices $i\in I$ are useful.
Thus
\begin{equation}\label{eqsizeI'}
|I'|\ge \tau|I|/8.
\end{equation}
Note that for every pair $i,i'\in I$
there exists an index $j'=j'(i,i')$ with $1\le j'<\ell$ and such that
$d_{G'}(Z_{j'}(i),Z_1(i'))=0$. (Otherwise we could take $i',i$ for~$i_1,i_2$.)
So in the graph $G'$ every vertex in $Z_1(i')$ has at most $(\ell-2+\xi/(1-\theta))L''$
neighbours in $Z_1(i)\cup\dots\cup Z_\ell(i)$.
Clearly, $2\le j'<\ell$ if $i$ is useful for~$i'$.

Given $i'\in I$, call another index $i\in I$ \emph{typical for $i'$}
if $d_{G'}(Z_j(i),Z_1(i'))>0$ for all $j\le \ell$ with $j\neq j'$.
Thus if $i$ is not typical for $i'$ then in the graph $G'$ every vertex
in $Z_1(i')$ has at most $(\ell-2)L''$ neighbours in $Z_1(i)\cup\dots\cup Z_\ell(i)$.
Given $i'\in I'$, we will now show that at least half of the $\ge \tau |I|/8$
indices $i$ which are useful for $i'$ are also typical.
Indeed, suppose not.
Consider any vertex $v\in Z_1(i')$ and look at its degree in~$G'$. We have that
\begin{align*}
d_{G'}(v) & \le |I|L''(1-\tau/16)(\ell-2+\xi/(1-\theta))+\frac{\tau}{16}|I|L''(\ell-2)
+\theta^{1/5}n\\
& \stackrel{(\ref{eqsizeI})}{\le} (\ell-2+\xi)k''L''-\tau \xi |I|L''/16+2\theta^{1/5}n
\le \delta(G)-\tau^2n<\delta(G'),
\end{align*}
a contradiction.
(Indeed, to see the first inequality use that the error bound $\theta^{1/5}n$
on the right hand side is a bound on the number of all those neighbours of~$v$
which lie in $V_0$ or in sets $Z_j(i)$ with $i\notin I$
(c.f.~(\ref{eqV0}) and~(\ref{eqsizeI})). Again, we have room to spare here.
To check the third inequality
use that~(\ref{sizeZji}) implies 
$n-|V_0|=(\ell-1+\xi/(1-\theta))L''k''\ge (\ell-1+\xi)L''k''$.
For the last inequality use that the Regularity lemma (Lemma~\ref{deg-reglemma})
implies $\delta(G')\ge \delta(G)-2dn$.)
This shows that for every $i'\in I'$ at least $\tau |I|/16$ indices $i\in I$
are both useful and typical for~$i'$.

Consider all the triples $i,i',j'$ such that $i\in I$, $i'\in I'$ and such that
$i$ is both useful and typical for~$i'$ and where $j'=j'(i,i')$ is
as defined after~(\ref{eqsizeI'}). 
It is easy to see that the number of such triples is at least
$\tau |I||I'|/16$. Thus there must be one pair $i,j'$ which occurs for at least
$\tau |I'|/(16\ell)$ indices $i'\in I'$. Let $I''$ denote the
set of all these indices~$i'$. So crudely
\begin{equation} \label{Iprimebound}
|I''| \stackrel{(\ref{eqsizeI'})}{\ge} \tau^3|I| \stackrel{(\ref{eqsizeI})}{\ge}
\tau^4 k''.
\end{equation}
Note that for each $i'\in I''$ there exists
a $j''$ such that $2\le j''\le\ell$ and $d_{G'}(Z_{j'}(i),Z_{j''}(i'))=0$.
(Otherwise we could take $i,i',j'$ for~$i_3,i_4,j^*$ since $i$ is both useful
and typical for~$i'$.)
So in the graph~$G'$ every vertex in $Z_{j'}(i)$ has at most $(\ell-2)L''$
neighbours in $Z_1(i')\cup\dots\cup Z_\ell(i')$.
 
Furthermore, for each $i''\in I\setminus (I''\cup \{i\})$ there exists
a $j'''$ such that $1\le j'''<\ell$ and $d_{G'}(Z_{j'}(i),Z_{j'''}(i''))=0$.
(Otherwise we could take $i,i',i'',j'$ for $i_5,i_6,i_7,j^\diamond$.)
Thus for each $i'' \in I\setminus (I''\cup \{i\})$ we can still say that in~$G'$
a vertex $v \in Z_{j'}(i)$ has at most $(\ell-2+\xi/(1-\theta))L''$ neighbours in
$Z_1(i'')\cup\dots\cup Z_\ell(i'')$.
Note that $v$ sends at most $\theta^{1/5}n$ edges to $V_0$ and to sets $Z_{i^*}$ with
$i^* \notin I$ by~(\ref{eqV0}) and~(\ref{eqsizeI}).
Together the above observations show that
\begin{align*}
d_{G'}(v) & \le (|I|-|I''|)(\ell-2+\xi/(1-\theta))L''+|I''|(\ell-2)L''+\theta^{1/5}n\\
& \le k''(\ell-2+\xi)L''-\xi|I''|L''/(1-\theta)+2\theta^{1/5}n
\stackrel{(\ref{Iprimebound})}{\le} \delta(G)-\tau^5 n<\delta(G'),
\end{align*}
a contradiction.

Thus we have shown that we can incorporate the first exceptional vertex $x$ by
removing at most three copies of $B^*$ which are as in~$(*)$.
Recall that this ensures that the subsets thus obtained from the $Z_j(i)$
satisfy conditions~(a) and~(c).
Next we proceed similarly with all other vertices in~$V_0$.
However, in order to ensure that in the end
condition~(b) is satisfied too, we need to be careful that we do
not remove to many vertices from a single set~$Z_j(i)$. So if the size of some
set $Z_j(i)$ becomes critical after we dealt with some vertex in~$V_0$,
then we exclude \emph{all} the sets $Z_1(i),\dots,Z_\ell(i)$ from consideration
when dealing with the remaining vertices in~$V_0$. 
The definition of the critical threshold in~(b) implies that we exclude at most
$$z|V_0|/(\theta^{1/7} L'')\stackrel{(\ref{eqV0})}{\le}
4z\sqrt{\theta}n/(\theta^{1/7}L'')\ll \theta^{1/3} k''
$$
indices $i$ in this way.
It is easy to check that this will not affect any of the above calculations
significantly. This completes the proof of the claim.

\medskip

\noindent
Recall that the sets $Z_j(i)$ were obtained by splitting the clusters
belonging to the copies $B_1,\dots, B_{k'}$ of $B^*$ in~$\B'$.
By taking out the copies of $B^*$ chosen in the above process we modified
these clusters.
For all $t\le k'$ and  all $j\le  \ell$ let $X_j(t)$ denote the union of all the
modified clusters belonging to the $j$th vertex class of~$B_t$.
Thus $X_j(t)=Y_j((t-1)z_1z+1)\cup\dots\cup Y_j(tz_1z)$. Moreover,~(\ref{sizeZji}) 
and (a)--(c) imply
that%
      \COMMENT{Indeed, (a) implies that all the bipartite graphs
between all the modified clusters $V,V'$ with
$V\subseteq X_j(t)$ and $V'\subseteq X_{j'}(t)$ are superregular.
But this in turn implies~(a$'$).}
\begin{itemize}
\item[(a$'$)] all the bipartite graphs $(X_j(t),X_{j'}(t))_{G'}$
are  $(\eps',d')$-superregular whenever $j\neq j'$,
\item[(b$'$)] $(1-\theta^{1/8})zL'\le |X_j(t)|\le (1-\theta)zL'$
for all $j<\ell$ and $(1-\theta^{1/8})z_1L'\le |X_\ell(t)|\le z_1L'$,
\item[(c$'$)] $|X_1(t)|=\dots=|X_{\ell-1}(t)|\le (1-\theta)|X_\ell(t)|/\xi$.
\end{itemize}

\subsection{Making the blow-ups of the $B\in \B'$ divisible by~$|B^*|$}\label{sec:div}

Given a subgraph $S\subseteq R$, we denote by $V_G(S)\subseteq V(G)$ the blow-up of~$V(S)$.
Thus $V_G(S)$ is the union of all the clusters which are vertices of~$S$. In particular,
$V_G(B_i)=X_1(i)\cup\dots\cup X_\ell(i)$. If $|V_G(B_i)|$ was divisible by~$|B^*|$
for each $B_i\in\B'$, then $\B'$ would correspond to a blown-up $B^*$
cover as required in the lemma. As already described at the end of
Section~\ref{sec:applyRG}, we will
achieve this by taking out a bounded number of further copies of $B^*$ from $G$.
For this, we define an auxiliary graph $F$ whose vertices are the
elements of $\B'$ and in which $B_i,B_j\in\B'$ are adjacent if the reduced
graph $R$ contains a copy of $K_\ell$ with one vertex in $B_i$ and $\ell-1$
vertices in $B_j$ or vice versa.

To motivate the definition of $F$, let us first consider the case when $F$ is
connected. If $B_i,B_j\in\B'$ are adjacent in $F$ then $G$ contains a copy of
$B^*$ with one
vertex in $V_G(B_i)$ and all the other vertices in $V_G(B_j)$ or vice versa.
In fact, we can even find $|B^*|-1$ disjoint such copies of $B^*$ in~$G$.
Taking out a suitable number of such copies (at most $|B^*|-1$), we can achieve that
the size of the subset of $V_G(B_i)$ obtained in this way is divisible by $|B^*|$.
Thus we can `shift the remainders mod~$|B^*|$' along a spanning tree of~$F$
to achieve that $|V_G(B)|$ is divisible by $|B^*|$ for each $B\in\B^*$.
(To see this, use that $\sum_{B\in\B'} |V_G(B)|$ is divisible by $|B^*|$
since $|G|$ is divisible by $|B^*|$.)

Let us next show that in the case when $\ell=2$ the graph $F$ is always
connected. If $\ell=2$, then $B_i,B_j\in \B'$ are joined in $F$ if
and only $R$ contains an edge between $B_i$ and~$B_j$.
Now suppose that $F$ is not connected and let $C$ be any component of~$F$.
Let $A\subseteq V(G)$ denote the union of all those clusters which
belong to some $B_i\in C$. Then in the current subgraph of $G'$ there are
no edges emanating from~$A$. As we have taken out at most $3|B^*||V_0|\le 12|B^*|\sqrt{\theta}n$
vertices in Section~\ref{sec:adjust} this implies that
$d_{G'}(A,V(G')\setminus A)\le \theta^{1/3}$.
Since $d_G(x)\le d_{G'}(x)+2dn$ for any $x\in G$ we have
$d_G(A,V(G)\setminus A)\le \theta^{1/3}+4d\ll \tau$,
a contradiction to condition~(ii) of Lemma~\ref{nonextremal}.

Thus in what follows we may assume that $\ell\ge 3$ and that $F$ is not connected. 
Let $\C$ denote the set of all components of~$F$.
Given a component $C$ of $F$, we denote by $V_R(C)\subseteq V(R)$ the set of all
those clusters which belong to some $B\in \B'$ with $B\in C$. Let $V_G(C)\subseteq V(G)$
denote the union of all clusters in $V_R(C)$.
We first show that we can take out a bounded number of copies of
$B^*$ from $G$ in order to make $|V_G(C)|$ divisible by $|B^*|$ for each $C\in\C$.
After that, we can `shift the remainders mod~$|B^*|$'
within each component $C\in\C$ along a spanning tree as indicated above to make
$|V_G(B)|$ divisible by~$|B^*|$ for each $B\in \B'$. For our argument, we will need
the following claim.

\medskip

\noindent
\textbf{Claim~1.} \emph{Let $C_1,C_2\in \C$ be distinct and let $a\in V_R(C_2)$. Then
$$|N_R(a)\cap V_R(C_1)|<
\frac{\ell-2}{\ell-1}|V_R(C_1)|.$$}

\smallskip

\noindent
Suppose not. Then there is some $B\in \B'$ such that $B\in C_1$ and
such that
$$
|N_R(a)\cap B|\ge \frac{\ell-2}{\ell-1} |B|> \frac{\ell-2}{\ell-1+\xi} |B|=(\ell-2)z.
$$
This implies
that $a$ has a neighbour in at least $\ell-1$ vertex classes of~$B$.
Thus $R$ contains a copy of $K_\ell$ which consists of $a$ together with $\ell-1$
of its neighbours in~$B$. But by definition of the auxiliary graph $F$, this means
that $B$ is adjacent in $F$ to the copy $B_i\in\B'$ that contains~$a$,
i.e.~$B$ and $B_i$ lie in the same component of $F$, a contradiction.
This completes the proof of Claim~1.

\medskip

\noindent
\textbf{Claim~2.} \emph{There exist a component $C'\in \C$, a copy $K$ of $K_\ell$ in $R$
and a vertex $a_0\in V(R)\setminus (V(K)\cup V_R(C'))$ such that $K$ meets $V_R(C')$ in
exactly one vertex and such that $a_0$ is joined to all the remaining vertices in $K$.
}

\smallskip

\noindent
As $\delta(R)> 1/2$,%
    \COMMENT{since we are in the case when $\ell\ge 3$}
there exists an edge $a_1a_2\in R$ which joins the vertex sets
corresponding to two different components of $F$, i.e.~there are distinct
$C_1,C_2\in\C$ such that $a_1\in V_R(C_1)$ and $a_2\in V_R(C_2)$.
Note that (\ref{eqmindegreduced}) implies that
\begin{equation}\label{eqweakdeltaR}
\delta(R)> \frac{\ell-2}{\ell-1}|R|.
\end{equation}
Thus the number of common neighbours of $a_1$ and $a_2$ in~$R$
is greater than
$\frac{\ell-3}{\ell-1}|R|.
$
To prove the claim, we will now distinguish two cases.

\medskip

\noindent
\textbf{Case~1.} \emph{More than
$\frac{\ell-3}{\ell-1}|V(R)\setminus V_R(C_1)|
$
common neighbours of $a_1$ and $a_2$ lie outside~$V_R(C_1)$.}

\smallskip

\noindent
Let $a_3$ be a common neighbour of $a_1$ and $a_2$ outside $V_R(C_1)$.
Claim~1 and (\ref{eqweakdeltaR})
together imply that the number of common neighbours of $a_1$, $a_2$ and $a_3$
outside $V_R(C_1)$ is more than
$$\frac{\ell-4}{\ell-1}|V(R)\setminus V_R(C_1)|.
$$
Choose such a common neighbour~$a_4$. Continuing in this way, we can obtain
distinct vertices $a_2,\dots,a_{\ell}$ outside $V_R(C_1)$ which together
with $a_1$ form a copy $K$ of $K_\ell$ in~$R$. As before, Claim~1
and (\ref{eqweakdeltaR}) together imply that
the number of common neighbours of $a_2,\dots,a_{\ell}$ outside
$V_R(C_1)$ is nonzero. Let
$a_0$ be such a common neighbour. Then Claim~2 holds with $C':=C_1$, $K$ and $a_0$.
Thus we may now consider

\medskip

\noindent
\textbf{Case~2.} \emph{More than
$\frac{\ell-3}{\ell-1}|V_R(C_1)|
$
common neighbours of $a_1$ and $a_2$ lie in~$V_R(C_1)$.}

\smallskip

\noindent
In this case we proceed similarly as in Case~1. However, this time we choose
$a_0,a_3,\dots,a_{\ell}$ inside $V_R(C_1)$. Indeed, this can be done since
Claim~1 and (\ref{eqweakdeltaR})
together imply  that each vertex in $V_R(C_1)$ has more than
$\frac{\ell-2}{\ell-1}|V_R(C_1)|
$
neighbours in~$V_R(C_1)$. Then Claim~2 holds with $C':=C_2$.

\medskip

\noindent
\textbf{Claim~3.} \emph{We can make $|V_G(B)|$ divisible by $|B^*|$ for all $B\in\B'$
by taking out at most $|\B'||B^*|$ disjoint copies of $B^*$ from~$G$.}

\smallskip

\noindent
We first take out some copies of $B^*$ from $G$ to achieve that $|V_G(C)|$ is
divisible by $|B^*|$ for each $C\in\C$. To do this we proceed as follows.
We apply Claim~2 to find a component $C_1\in\C$,
a copy $K$ of $K_\ell$ in $R$ and a vertex $a_0\in V(R)\setminus (V(K)\cup V_R(C_1))$
such that $K$ meets $V_R(C_1)$ in exactly one vertex, $a_1$ say, and such that
$a_0$ is joined to all vertices in $K-a_1$. Thus $G$ contains a copy $B'$ of $B^*$ which 
has exactly 
one vertex $x\in V_G(C_1)$ and whose other vertices lie in clusters belonging to
$V(K-a_1)\cup \{a_0\}$. (Indeed, we can choose the vertices of $B'$ lying in the same vertex
class as $x$ in the cluster $a_0$ and the vertices lying in other vertex
classes in the clusters belonging to $K-a_1$.) In fact, $G$ contains $|B^*|-1$ (say)
disjoint such copies of~$B^*$. Now suppose that
$|V_G(C_1)|\equiv j\mod |B^*|$. Then we take out $j$ disjoint such copies of $B^*$ from $G$
to achieve that $|V_G(C_1)|$ is divisible by $|B^*|$. Next we consider the graphs
$F_1:=F-V(C_1)$ and $R_1:=R-V_R(C_1)$ instead of $F$ and $R$.
Claim~1 and (\ref{eqweakdeltaR}) together imply that
$\delta(R_1)>\frac{\ell-2}{\ell-1}|R_1|.
$
Now suppose that $|\C|\ge 3$. Then similarly as in the proof of Claim~2 one can find
a component $C_2\in\C\setminus\{C_1\}$,
a copy $K'$ of $K_\ell$ in $R_1$ and a vertex $a'_0\in V(R_1)\setminus (V(K')\cup V_R(C_2))$
such that $K'$ meets $V_R(C_2)$ in exactly one vertex, $a_2$ say, and such that
$a'_0$ is joined to all vertices in $K-a_2$. As before, we take out at most $|B^*|-1$ copies of $B^*$
from $G$ to achieve that $|V_G(C_2)|$ is divisible by~$|B^*|$. As $|G|$ was divisible by
$|B^*|$, we can continue in this fashion to achieve that $|V_G(C)|$ is divisible
by $|B^*|$ for all components $C\in\C$. In this process, we have to take out at most
$(|\C|-1)(|B^*|-1)$ copies of $B^*$ from $G$. Now we consider each component
$C\in \C$ separately. By proceeding as in the connected case for each 
$C$ and taking out at most $(|C|-1)(|B^*|-1)$ further copies of
$B^*$ from $G$ in each case, 
we can make $|V_G(B)|$ divisible by $|B^*|$ for each $B\in \B'$.
Hence, in total, we have taken out at most $(|\C|-1)(|B^*|-1)+(|\B'|-|\C|)(|B^*|-1)\le
|\B'||B^*|$ copies of $B^*$ from~$G$.

$\B'$ now corresponds to a blown-up $B^*$-cover as desired in the lemma,
except that it has $k'\le k_1$ elements. (Recall that $k_1$ was defined
in~(\ref{eqdefk1}).) But by considering random partitions,
it is easy to see that one can split these elements to obtain a blown-up $B^*$-cover
as required. 
The $B^*$-packing $\B^*$ in Lemma~\ref{nonextremal}
consists of all the copies of $B^*$ taken out during the proof.
Thus $|\B^*|\le 3|V_0|+|\B'||B^*|\le 12\sqrt{\theta}n+|\B'||B^*|\le \theta^{1/3}n$,
as desired.  

\section{Packings in complete $\ell$-partite graphs}\label{sec:complete}

In this section, we prove several results which together imply Lemma~\ref{corcomplete1}.
However, almost all of the results of this section are also used directly in 
Section~\ref{sec:extremal}.

Clearly, a complete $\ell$-partite graph has a perfect $H$-packing if all its
vertex classes have equal size which is divisible by~$|H|$.
Together the following two lemmas show that if $\hcf(H)=1$ then we still
have a perfect $H$-packing
if the sizes of the vertex classes are permitted to deviate slightly. 
(By Proposition~\ref{bestposs1} this is false if $\chi(H)\ge 3$ and
$\hcf(H)\neq 1$ or if $\chi(H)=2$ and $\hcf_\chi(H)>2$.)
In Lemma~\ref{completemove} we first consider the case when $\hcf_\chi(H)=1$.
In Lemma~\ref{completebipmove2} we then deal with the remaining case
(i.e.~when $H$ is bipartite, $\hcf(H)=1$ but $\hcf_\chi(H)=2$).

\begin{lemma}\label{completemove}
Suppose that $H$ is a graph of chromatic number $\ell\ge 2$ such that
$\hcf_\chi(H)=1$. Let $B^*$ be the bottlegraph
assigned to~$H$. Let $D'\gg |H|$ be an integer divisible by~$|H|$.
Let $a$ be an integer such that $|a|\le |B^*|$.
Given $1\le i_1<i_2\le \ell$, let $G$ be a complete $\ell$-partite graph
with vertex classes $U_1,\dots,U_\ell$ such that $|U_{i_1}|=D'+a$, $|U_{i_2}|=D'-a$
and $|U_r|=D'$ for all $r\neq i_1,i_2$. Then $G$ contains a perfect
$H$-packing.
\end{lemma}
\proof
Let $q$ denote the number of optimal colourings of~$H$.
First note that by taking out at most $q\ell!$ disjoint copies of $H$
from $G$ we may assume that $D'$ is divisible by~$q(\ell-1)!|H|$.%
     \COMMENT{Indeed, let $0\le t<q(\ell-1)!$ be an integer such that
$D'=|H|(q(\ell-1)!k+t)$ where $k\in\mathbb{N}$. Take out copies of $H$ in $t$
steps. In each step we take out $\ell$ copies such that from each $U_i$
we remove $|H|$ vertices (in each step).}
Consider the complete $\ell$-partite graph $G'$ whose vertex classes
$U'_1,\dots,U'_\ell$ all have size~$D'$. Thus $|G|=|G'|$.
We will think of $G$ and $G'$ as two graphs
on the same vertex set whose vertex classes are roughly identical.
Recall that $G'$ has a perfect $H$-packing.

Our aim is to
choose a suitable such $H$-packing $\cH'$ in $G'$ and to show that it can be modified
into a perfect $H$-packing of $G$. To choose
$\cH'$, we will consider all optimal colourings
of~$H$. Let $c^1,\dots,c^q$ be all these colourings.
Let $x^j_1\le x^j_2\le \dots\le x^j_{\ell}$
denote the sizes of the colour classes of~$c^j$. 
Put
$$
k:=\frac{D'}{q(\ell-1)!|H|}.
$$
Let $S_\ell$ denote the set of all permutations of $\{1,\dots,\ell\}$.
Given $j\le q$, let $\cH'_j$ be an $H$-packing in $G'$ which, for all
$s\in S_\ell$, contains precisely $k$ copies of $H$ which in the colouring $c^j$
have their $s(i)$th colour class in $U'_i$ (for all $i=1,\dots,\ell$).
Thus $\cH'_j$ consists of $\ell!k$ copies of $H$ and covers precisely
$k(\ell-1)!|H|$ vertices in each vertex class $U'_i$ of~$G'$.
Moreover, we choose all the $\cH'_j$ to be disjoint from each other.
Thus the union $\cH'$ of $\cH'_1,\dots,\cH'_q$ is a perfect
$H$-packing in~$G'$.

We will now show that $\cH'$ can be modified into a perfect $H$-packing of~$G$.
Roughly, the reason why this can be done is the following. Clearly,
we may assume that $a>0$.
So $\cH'$ has less than $|U_{i_1}|$ vertices in its $i_1$th vertex class and more than
$|U_{i_2}|$ vertices in its $i_2$th vertex class. We will modify $\cH'$ slightly by
interchanging some vertex classes in some copies of $H$ in $\cH'$. As $\hcf_\chi(H)=1$
this can be done in such a way that the $H$-packing obtained from $\cH'$ covers one
vertex more in its $i_1$th vertex class than $\cH'$ and one vertex less in its $i_2$th
vertex class. Continuing in this fashion we obtain an $H$-packing which covers the
correct number of vertices in each vertex class.

For all $j\le q$ and $r<\ell$ put $d^j_r:=x^j_{r+1}-x^j_r$.
Since $\hcf_\chi(H)=1$, we can find $b^j_r\in\mathbb{Z}$ such that
$$
1=\sum_{j=1}^q\sum_{r=1}^{\ell-1} b^j_rd^j_r.
$$
(Here we take $b^j_r:=0$ if $d^j_r=0$.)
In order to modify the $H$-packing $\cH'$ we proceed as follows.
For all $j\le q$ and $r<\ell$ we consider $b^j_r$.
If $b^j_r\ge 0$ we choose $b^j_r$ of the copies of $H$ in $\cH'_j\subseteq \cH'$
which in the colouring $c^j$ have their $r$th vertex class in $U'_{i_1}$
and their $(r+1)$th vertex class in $U'_{i_2}$. We change each of these
copies of $H$ such that they now have their $r$th vertex class in $U'_{i_2}$ and
its $(r+1)$th vertex class in $U'_{i_1}$. All the other vertices remain
unchanged. Note the number of vertices in the $i_1$th vertex class covered by
this new $H$-packing increases by $b^j_rd^j_r$ whereas the number of covered vertices
in the $i_2$th vertex class decreases by~$b^j_rd^j_r$. 

If $b^j_r< 0$ we choose $|b^j_r|$ of the copies of $H$ in $\cH'_j$
which in the colouring $c^j$ have their $r$th vertex class in $U'_{i_2}$
and their $(r+1)$th vertex class in $U'_{i_1}$. This time, we change each of these
copies of~$H$ such that they now have their $r$th vertex class in $U'_{i_1}$ and
its $(r+1)$th vertex class in~$U'_{i_2}$.
Note that all these copies of $H$ will automatically be distinct for different
pairs $j,r$. 

Let $\cH^*$ denote the modified $H$-packing obtained by proceeding as
described above $a$ times (where all the copies of $H$ which we change are
chosen to be distinct).%
     \COMMENT{Indeed, this can be done since $D'\gg |H|$, $|a|\le |B^*|$ and thus
$k\gg |b^j_r| a$ for all $j,r$.}
We have to check that
$\cH^*$ is a perfect $H$-packing of~$G$. For all $i\le \ell$ let $n_i$ denote
the number of vertices in the $i$th vertex class covered by~$\cH^*$. 
Thus $n_i=D'$ whenever $i\neq i_1,i_2$.
We have to check that $n_{i_1}=D'+a$ and $n_{i_2}=D'-a$. But
$$
n_{i_1}=D'+a\sum_{j=1}^q\sum_{r=1}^{\ell-1} b^j_rd^j_r=D'+a
$$
and
$$
n_{i_2}=D'-a\sum_{j=1}^q\sum_{r=1}^{\ell-1} b^j_rd^j_r=D'-a,
$$ 
as required.
\endproof

\begin{lemma}\label{completebipmove2}
Suppose that $H$ is a bipartite graph such that $\hcf_c(H)=1$ and
$\hcf_\chi(H)=2$. Let $B^*$ be the bottlegraph
assigned to~$H$. Let $D'\gg |H|$ be an integer divisible by~$|H|$.
Let $a$ be an integer such that $|a|\le |B^*|$.
Let $G$ be a complete bipartite graph
with vertex classes $U_1$ and $U_2$ such that $|U_1|=D'+a$ and $|U_2|=D'-a$.
Then $G$ contains a perfect $H$-packing.
\end{lemma}
\proof
Clearly, we may assume that $a>0$. Our first aim is to take out a small number
of disjoint copies of $H$ from $G$ to obtain sets $U'_i\subseteq U_i$ with
$|U'_1|=|U'_2|$. To do this, we will
use the fact that $\hcf_\chi(H)=2$. So let $c^1,\dots,c^q$ be all the optimal
colourings of~$H$. Let $x^j_1\le x^j_2$ denote the sizes of the colour
classes of~$c^j$. Since $\hcf_\chi(H)=2$, we can find $b^j\in\mathbb{Z}$
such that $2=\sum_{j=1}^q b^j (x^j_2-x^j_1)$. (Here we take $b^j:=0$ if $x^j_2=x^j_1$.)
For each $j=1,\dots,q$ in turn we take out $a|b^j|$ copies of
$H$ from~$G$. If $b^j\ge 0$ each of these $a|b^j|$ copies will meet
$U_1$ in $x^j_2$ vertices and $U_2$ in $x^j_1$ vertices. If $b^j<0$
then each of these copies will have $x^j_1$ vertices in
$U_1$ and $x^j_2$ vertices in~$U_2$. We choose all these copies of $H$
to be disjoint. It is easy to check that the subsets $U'_1$ and $U'_2$
obtained from $U_1$ and $U_2$ in this way have the same size,
$u'$ say.%
     \COMMENT{Indeed, $|U'_1|-|U'_2|=2a-a\sum_j b^j(x^j_2-x^j_1)=0$.}
Note that $|H|$ divides $2u'$ since $|H|$ divides $|U_1|+|U_2|$.
Also observe that if $|H|$ even divides $u'$,
then $G[U'_1\cup U'_2]$ (and thus also~$G$ itself) has a perfect $H$-packing.
So we may assume that $|H|$ does not divide~$u'$ but that it does divide $2u'$.
Hence $|H|$ is even and $u'=|H|k/2$ where $k$ is an odd integer.

We will now use the fact that $\hcf_c(H)=1$ to show that we can take out further copies
of $H$ from $G$ to achieve that the subsets $U''_1$ and $U''_2$ obtained in this way
have the same size and that this size is divisible by~$|H|$.
As $\hcf_c(H)=1$ there exists
a component $C$ of $H$ such that $|C|$ is odd. Using the fact that $|H|$ is even
and thus $|H-C|$ is odd it is easy to see that there exists a $2$-colouring
of $H$ whose colour classes both have odd size and another $2$-colouring
whose colour classes both have even size.%
      \COMMENT{Indeed, let $a_1$ and $a_2$ denote the
sizes of the colour classes of $C$ (in some 2-colouring). We may assume that $a_1$ is
odd and $a_2$ is even. Since $|H|$ is even it follows that $|H-C|$ is odd.
Let $a_3$ and $a_4$ denote the sizes of the colour classes of $H-C$ (in some $2$-colouring)
such that $a_3$ is odd and $a_4$ is even. Thus there exists a 2-colouring of
$H$ whose colour classes have size $a_1+a_3$ (even) and $a_2+a_4$ (even).
Another colouring has colour classes of size $a_1+a_4$ (odd) and
$a_2+a_3$ (odd).}
We may assume that $c^1$ and $c^2$ are such colourings, i.e. that both $x^1_1$ and
$x^1_2$ are odd and both $x^2_1$ and $x^2_2$ are even. 
Let $k_1:=|H|/2-x^2_1$ and $k_2:=|H|/2-x^1_1$. Take out $k_1$ copies of
$H$ with $x^1_1$ vertices in $U'_1$ and $x^1_2$ vertices in~$U'_2$.
Then take out $k_2$ copies of $H$ with $x^2_2$ vertices in $U'_1$ and $x^2_1$
vertices in~$U'_2$. Let $U''_1$ and $U''_2$ denote the subsets obtained
from $U'_1$ and $U'_2$ in this way. It is easy to check that%
     \COMMENT{Indeed,
\begin{align*}
|U''_1|&= u'-(\frac{|H|}{2}-x^2_1)x^1_1-(\frac{|H|}{2}-x^1_1)x^2_2=
u'-(\frac{|H|}{2}-x^2_1)x^1_1-(\frac{|H|}{2}-x^1_1)(|H|-x^2_1)\\
& =u'-\frac{|H|}{2}(x^1_1-x^2_1+|H|-2x^1_1)
\end{align*}
and
\begin{align*}
|U''_2|&= u'-(\frac{|H|}{2}-x^2_1)x^1_2-(\frac{|H|}{2}-x^1_1)x^2_1=
u'-(\frac{|H|}{2}-x^2_1)(|H|-x^1_1)-(\frac{|H|}{2}-x^1_1)x^2_1\\
& =u'-\frac{|H|}{2}(-x^1_1-2x^2_1+|H|+x^2_1).
\end{align*}}
$$|U''_1|=|U''_2|=u'-\frac{|H|}{2}(|H|-x^1_1-x^2_1)=\frac{|H|}{2}(k-|H|+x^1_1+x^2_1).$$
But $k-|H|+x^1_1+x^2_1$ is even and so $|U''_1|$ is divisible by $|H|$, as desired.
\endproof

The next lemma is an analogue of Lemma~\ref{completebipmove2}
for perfect $H$-packings in graphs which are the disjoint union of two cliques.
It will be needed in the proof of Lemma~\ref{extremal2}.

\begin{lemma}\label{completebipmove}
Suppose that $H$ is a bipartite graph such that
$\hcf_c(H)=1$. Let $B^*$ be the bottlegraph
assigned to~$H$. Let $D'\gg |H|$ be an integer divisible by~$|H|$.
Let $a$ be an integer such that $|a|\le |B^*|$.
Let $G$ be the disjoint union of two cliques 
of order $D'+a$ and  $D'-a$ respectively. Then $G$ contains a perfect
$H$-packing.
\end{lemma}
\proof
Clearly, we may assume that $a>0$. Let $G_1$ be the clique of order $D'+a$ and
let $G_2$ be the clique of order $D'-a$. Our aim is to
take out disjoint copies of $H$
from $G$ in order to obtain subcliques $G'_i\subseteq G_i$ such that
$|G'_i|$ is divisible by $|H|$ for both $i=1,2$.
(Then each $G'_i$ (and thus also~$G$ itself) has a perfect $H$-packing.)

Let $C_1 \le \dots \le C_s$ denote the components of $H$.
Since $\hcf_c(H)=1$, we can find $b^j\in\mathbb{Z}$ such that
$
1=\sum_{j=1}^s b^j |C_j|.
$
For each $j=1,\dots,s$ in turn we take out $a|b^j|$ copies of
$H$ from~$G$. If $b^j\ge 0$ each of these $a|b^j|$ copies will meet
$G_1$ in $C_j$ and $G_2$ in~$H-C_j$. If $b^j<0$
then each of these copies will meet $G_1$ in $H-C_j$ and $G_2$ in~$C_j$.
We choose all these copies of $H$ to be disjoint. Let $G'_1$ and $G'_2$
denote the subcliques obtained from $G_1$ and $G_2$ in this way.
Then
\begin{align*}
|G'_1|=D'+a-a\sum_{j:\, b^j\ge 0} b_j|C_j|+a \sum_{j:\, b^j< 0} b^j(|H|-|C_j|)
= D'+a \sum_{j:\, b^j< 0} b^j|H|.
\end{align*}
Thus $|H|$ divides $|G'_1|$. Similarly one can show that $|H|$ divides $|G'_2|$.
\endproof
 
The following lemma states that if $G$ is a complete $\ell$-partite graph which is 
very close to being bottle-shaped then $G$ contains a perfect $H$-packing
as long as the ratio of the smallest to the largest vertex class a bit larger than
in the bottlegraph $B^*$ of $H$. (The terms involving $D'$ in (i) and (ii)
ensure that the latter condition holds.) 
\begin{lemma}\label{completeconst}
Suppose that $H$ is a graph of chromatic number $\ell\ge 2$ such that
$\hcf(H)=1$. Let $B^*$ be the bottlegraph assigned to~$H$.
Let $z$ and $z_1$ be as defined in~$(\ref{eqdefxi})$.
Let $D'\gg |H|$ be an integer divisible by~$|B^*|$.
Let $G$ be a complete $\ell$-partite graph with vertex classes $U_1,\dots,U_\ell$
whose order $n\gg D'$ is divisible by~$|B^*|$. Let $u_i:=|U_i|$ for all~$i$.
Suppose that 
\begin{itemize}
\item[(i)] $|\, (u_i-D')-z(n-\ell D')/|B^*|\,|\le |B^*|$ for all $i<\ell$ and
\item[(ii)] $|\, (u_\ell-D')-z_1(n-\ell D')/|B^*|\,|\le |B^*|$.
\end{itemize}
Then one can take out $\ell D'/|H|$ disjoint copies of $H$ from $G$ to obtain
a subgraph $G^*\subseteq G$ such that, writing $n^*:=|G^*|=n-\ell D'$ and
$u^*_i:=|U_i\cap V(G^*)|$, we have $u^*_i=zn^*/|B^*|$ for all $i<\ell$
and $u^*_\ell=z_1n^*/|B^*|$. So in particular, $G^*$ contains a perfect
$B^*$-packing and thus $G$ contains a perfect $H$-packing.
\end{lemma}
\proof
Let us first consider the case when $\hcf_\chi(H)=1$.
Put $a_i:=(u_i-D')-z(n-\ell D')/|B^*|$ for all $i<\ell$ and
$a_\ell:=(u_\ell-D')-z_1(n-\ell D')/|B^*|$. Thus $\sum_{i=1}^\ell a_i=0$
and $|a_i|\le |B^*|$ for all~$i\le \ell$. Consider the complete
$\ell$-partite graph $G'$ whose $i$th vertex class has size $D'+a_i$.
By repeated applications of Lemma~\ref{completemove} one can show that
this graph has a perfect $H$-packing~$\cH$. View $G'$ as a subgraph of $G$
such that the $i$th vertex class of $G'$ lies in~$U_i$.
Then the subgraph $G^*$ obtained from $G$ by removing all the copies of
$H$ in $\cH$ (and thus deleting precisely the vertices in $V(G')\subseteq V(G)$)
is as required in the lemma.%
     \COMMENT{Note that $n^*=n-\ell D'$ is divisible by $|B^*|$
since both $n$ and $D'$ are divisible by~$|B^*|$. Thus $G^*$ is a `multiple' of $B^*$
and hence contains a perfect $B^*$-packing.}
In the remaining case when $\hcf_\chi(H)\neq 1$
(and thus $\chi(H)=2$, $\hcf_c(H)=1$ and $\hcf_\chi(H)=2$) we proceed similarly
except that we now apply Lemma~\ref{completebipmove2} instead of
Lemma~\ref{completemove}.
\endproof
The next lemma shows that we can achieve the conditions in the setup of
Lemma~\ref{completeconst} when larger deviations from the bottle-shape are allowed.
\begin{lemma}\label{completeapprox}
Let $H$, $B^*$, $z$, $z_1$, $G$, $U_i$, $u_i$ and $D'$ be defined as in the previous lemma,
except that we now no longer assume that $\hcf(H)=1$ and that $G$ satisfies~(i) and~(ii)
and we only require
$D'\ge 0$ to be any integer divisible by~$|B^*|$. Suppose that $\chi_{cr}(H)<\ell$ and
$a_i:=z(n-\ell D')/|B^*|-(u_i-D')\ge 0$ where $a_i\le n/(\ell^3|B^*|^2)$
for all $i<\ell$. Then one can take out at most $\ell^2 \sum_{i=1}^{\ell-1} a_i$
disjoint copies of $B^*$ from $G$ to obtain
a subgraph $G^*\subseteq G$ such that, writing $n^*:=|G^*|$ and
$u^*_i:=|U_i\cap V(G^*)|$, 
conditions (i) and (ii) of Lemma~\ref{completeconst} hold with $n$ replaced 
by $n^*$ and with $u_i$ replaced by~$u^*_i$.
\end{lemma}
\proof
First note that we only need to consider the case when $D'=0$. Indeed,
to reduce the general case, suppose that Lemma~\ref{completeapprox} holds
if $D'=0$. Now instead of~$G$ we consider the graph $G'$ obtained from
$G$ by removing $D'$ vertices from each vertex class. Apply Lemma~\ref{completeapprox} to $G'$
to obtain a graph $G^*\subseteq G'$. Let $U^*_1,\dots, U^*_\ell$ denote
the vertex classes of~$G^*$. Then the vertex classes obtained from
the $U^*_i$ by adding the $D'$ vertices span a subgraph of $G$ as desired in
the lemma. We may also assume that $u_1\le \dots\le u_{\ell-1}$.
Let $k:=n/|B^*|$ and let $\xi$ be as defined in~(\ref{eqdefxi}). Note that
\begin{equation}\label{eqsizeuell}
\quad u_i=kz-a_i  \mbox{ for } i<\ell \quad \mbox{ and } \quad u_\ell=\xi zk+\sum_{i=1}^{\ell-1} a_i .
\end{equation}
We now take out disjoint copies of $B^*$ from $G$ in order to achieve that
the subsets of $U_1,\dots,U_{\ell-1}$
thus obtained have almost the same size. More precisely, we proceed as follows.
For every $i=1,\dots,\ell-2$ let $r_i:=\lfloor(u_{\ell-1}-u_i)/(z-z_1)\rfloor$.
Put
\begin{align*}
r:=\sum_{i=1}^{\ell-2} r_i
\le \frac{(\ell-2)(kz-a_{\ell-1})-\sum_{i=1}^{\ell-2} (kz-a_i)}{z-z_1}
=\frac{\sum_{i=1}^{\ell-2} a_i-(\ell-2)a_{\ell-1}}{z-z_1}.
\end{align*}
For every $i=1,\dots,\ell-2$ in turn remove $r_i$ copies of $B^*$
from $G$, each having $z_1$ vertices in $U_i$ and $z$ vertices in every other
set~$U_j$. Then the subsets $U'_i$ obtained from the $U_i$ in this way
satisfy $0\le |U'_{\ell-1}|-|U'_i|<z-z_1$ for all $i=1,\dots,\ell-2$ and
$$
|U'_\ell|-\xi |U'_{\ell-1}|=u_\ell-\xi u_{\ell-1}-r(z-z_1)
\stackrel{(\ref{eqsizeuell})}{\ge}(\ell-1+\xi)a_{\ell-1}\ge 0.
$$
Next we will take out further copies of $B^*$ from $G$ in order to achieve that the
size of the $\ell$th vertex class is about $\xi$-times as large as the size of any
other vertex class. In each step we remove $\ell-1$
copies of $B^*$, for every $i=1,\dots,\ell-1$ one copy having $z_1$ vertices in
the $i$th vertex class and $z$ vertices in each other class. A straightforward
calculation shows that after
$$
\left\lfloor \frac{|U'_\ell|-\xi |U'_{\ell-1}|}{(\ell-1)z-(\ell-2)z_1-\xi z_1} \right\rfloor
$$
steps the subsets $U^*_i$ obtained from the $U'_i$ in this way 
span a subgraph as required in the lemma.%
     \COMMENT{First note that all the $U^*_i$ are non-empty. Indeed,
in total we took out at most
\begin{align*}
r+(\ell-1)\frac{|U'_\ell|-\xi |U'_{\ell-1}|}{(\ell-1)z-(\ell-2)z_1-\xi z_1}
\le \sum_{i=1}^{\ell-2} a_i+\ell^2 a_{\ell-1}\le \ell^2\sum_{i=1}^{\ell-1} a_i
\le \frac{n}{|B^*|^2}
\end{align*}
copies of $H$. Here we also used that the lower bound on $|U'_\ell|-\xi |U'_{\ell-1}|$
is almost an equality. Thus from each $U_i$ we deleted at most $zn/|B^*|^2\le zn/2|B^*|\le 2u_i/3$
vertices. To check that the size of $U^*_\ell$ is ok let
$x:=\lfloor (|U'_\ell|-\xi |U'_{\ell-1}|)/((\ell-1)z-(\ell-2)z_1-\xi z_1)\rfloor$
and note that
\begin{align*}
u^*_\ell-\xi u^*_{\ell-1} & =u'_\ell-x(\ell-1)z-\xi (u'_{\ell-1}-x((\ell-2)z+z_1))\\
& = u'_\ell-\xi u'_{\ell-1}-x[(\ell-1)z-(\ell-2)z_1-\xi z_1].
\end{align*}
Thus $0\le u^*_\ell-\xi u^*_{\ell-1}\le (\ell-1)z$ and $0\le u^*_{\ell-1}-u^*_i\le z$
for all $i<\ell-1$. Let us now show that this implies that the $u^*_i$ are as
required in the lemma. Indeed, we have that
$(\ell-2)z\ge (\ell-2)u^*_{\ell-1}-\sum_{i=1}^{\ell-2}u^*_i\ge 0$ and
thus
$$(\ell-2)z\ge (\ell-2+\xi)u^*_{\ell-1}-\sum_{i=1}^{\ell-2}u^*_i-u^*_\ell
\ge-(\ell-1)z.$$
Note that the term in the middle equals
$(\ell-1+\xi)u^*_{\ell-1}-n^*.$
Thus $-|B^*|\le -\frac{\ell-2}{(\ell-1+\xi)}z\le \frac{zn^*}{|B^*|}-u^*_{\ell-1}=
\frac{n^*}{\ell-1+\xi}-u^*_{\ell-1}\le \frac{\ell-1}{\ell-1+\xi}z\le |B^*|.$
Similarly, we get that
\begin{align*}
u^*_\ell & \le  \xi u^*_{\ell-1}+(\ell-1)z
\le \xi\frac{zn^*}{|B^*|}+\frac{(\ell-2)\xi z}{\ell-1+\xi}+(\ell-1)z
\le\frac{\xi n^*}{\ell-1+\xi}+|B^*|
\end{align*}
and so
$$u^*_\ell-\frac{\xi n^*}{\ell-1+\xi}\le |B^*|.$$
Finally, we have that $u^*_\ell\ge \xi u^*_{\ell-1}\ge \frac{\xi n^*}{\ell-1+\xi}-\xi|B^*|$
and so $u^*_\ell-\frac{\xi n^*}{\ell-1+\xi}\ge -|B^*|$.
To get the bounds for the other $u_i^*$, note the above calculations have some room to 
spare}
\endproof

\removelastskip\penalty55\medskip\noindent{\bf Proof of Lemma~\ref{corcomplete1}. }
Let $D'$ be an integer as in Lemma~\ref{completeconst}.
Consider the graph $F$ given in Lemma~\ref{corcomplete1}. By taking out
at most $\ell-2$ disjoint copies of $H$ from $F$ if necessary, we may
assume that $|F|$ is divisible by~$|B^*|$. It is easy to
check that%
     \COMMENT{Indeed, we may assume that the vertex classes in
Lemma~\ref{corcomplete1} satisfy $(1-d)u_{\ell-1}\le u_i\le u_{\ell-1}$
for all $i<\ell-1$. Together with the fact that
$(1-\beta^{1/10})u_\ell\le \xi u_i\le (1-\beta)u_\ell$
this implies that
$$
 \frac{zn}{|B^*|}-u_{\ell-1} = \frac{1}{\ell-1+\xi}(\sum_{i=1}^{\ell} u_i)-u_{\ell-1}
\le \frac{\ell-1}{\ell-1+\xi} u_{\ell-1} + \frac{\xi}{\ell-1+\xi} 
\frac{u_{\ell-1}}{1-\beta^{1/10}} -u_{\ell-1}
\le \beta^{1/20} u_{ \ell-1} \le \beta^{1/20}n
$$
The lower bound follows in the same way, just using the other side of the previous inequalities.
Including the $D'$ in the inequalities doesn't change things significantly
if $n\gg D'$.
}
$F$ satisfies the conditions in Lemma~\ref{completeapprox} if $|F|\gg D'$.
Thus Lemmas~\ref{completeapprox} and~\ref{completeconst} together
imply Lemma~\ref{corcomplete1}.
\endproof

\section{Proof of the extremal cases}\label{sec:extremal}
In most of the extremal cases, we know that $G$ contains several large almost independent
sets $A_1,\dots,A_q$ where $1\le q <\ell$. In the preliminary Lemma~\ref{exceptionalvs}
we show that we can modify
the $A_i$ slightly to obtain sets $A_1^*,\dots,A_{q}^*$ which together with
$V(G)\setminus \bigcup_{i=1}^q A_i^*$ induce an almost complete
$(q+1)$-partite graph.

In the proof of Lemma~\ref{exceptionalvs} below we need the following observation.

\begin{lemma}\label{disjointstars}
Let $i$ be a positive integer and let $G$ be a graph of order $n$ whose
average degree satisfies $d:=d(G)\ge 2i$. Then $G$ contains at least
$$\frac{dn}{4(i+1)\Delta(G)}$$
disjoint $i$-stars.
\end{lemma}
\proof
Let $k:=\lceil dn/(4(i+1)\Delta(G))\rceil$. We take out the disjoint $i$-stars greedily.
So in each step we delete $i+1$ vertices and thus at most $(i+1)\Delta(G)$ edges.
So after $<k$ steps the remaining subgraph $G'$ of $G$ has at least
$e(G)-k(i+1)\Delta(G)\ge dn/4$ edges and thus $d(G')\ge d/2\ge i$. So
$G'$ contains an $i$-star. This shows that the number of disjoint $i$-stars
we can find greedily is at least~$k$.
\endproof

\begin{lemma}\label{exceptionalvs}
Suppose that $H$ is a graph of chromatic number $\ell\ge 2$
such that $\chi_{cr}(H)<\ell$.
Let $B^*$ denote the bottlegraph assigned to~$H$.
Let $\xi$, $z$ and $z_1$ be as defined in~$(\ref{eqdefxi})$ and 
let $0<\tau \ll \xi, 1-\xi, 1/|B^*|$.
Let $|B^*|\ll D'\ll C$ be integers such that $D'$ is divisible by $|B^*|$.
Let $G$ be a graph whose order $n\gg C,1/\tau$ is divisible by $|B^*|$
and whose minimum degree satisfies
$\delta(G)\ge (1-\frac{1}{\chi_{cr}(H)})n+C$. Furthermore, suppose
that for some $1\le q< \ell$ there are disjoint sets $A_1,\dots,A_q\subseteq V(G)$
which satisfy $|A_i|=(n-2\ell D')z/|B^*|+2D'$ and $d(A_i)\le \tau$. Let 
$A_{q+1}:=V(G)\setminus \bigcup_{i=1}^q A_i$. Then
there are sets $A^*_1,\dots,A^*_{q+1}$ which satisfy the following properties:
\begin{itemize}
\item[{\rm (i)}] Let $A^*$ denote the union of $A^*_1,\dots,A^*_{q+1}$ and put
$n^*:=|A^*|$. Then $G-A^*$ has a perfect $H$-packing.
Moreover $n-n^*\le \tau^{3/5} n$.
\item[{\rm (ii)}] $|A^*_i|=(n^*-\ell D')z/|B^*|+D'$ and $d(A_i^*)\le \tau^{2/5}$
for all $i\le q$.
\item[{\rm (iii)}] For all $i,j\le q+1$ with $j\neq i$ each vertex in $A^*_i$
has at least $(1-\tau^{1/5})|A^*_j|$ neighbours in $A^*_j$.
\end{itemize}
\end{lemma}
\proof
First note that the minimum degree condition on $G$ and~(\ref{eqchicr}) 
imply that the neighbourhood
of each vertex $x\in G$ can avoid almost $|A_1|=\dots=|A_q|$ vertices of $G$ but no more.
Given an index $i\le q+1$, we call a vertex $x\in A_i$ \emph{$i$-bad} if
$x$ has at least $\tau^{1/3}|A_i|$ neighbours in $A_i$. Since $d(A_i)\le \tau$
for each $i\le q$, for such $i$'s the number of $i$-bad vertices is at most
$\tau^{2/3}|A_i|$. Call a vertex $x\in A_i$ \emph{$i$-useless} if $x$ has at most
$(1-\tau^{1/4})|A_j|$ neighbours in~$A_j$ for some $j\neq i$. Thus if $i\le q$
every $i$-useless vertex is also $i$-bad.%
     \COMMENT{Indeed to see this (and also for later arguments) it is helpful to
note that
$$|A_1|=\frac{(n-2\ell D')z}{(\ell-1+\xi)z}+2D'=
\frac{n}{\ell-1+\xi}-\frac{2D(1-\xi)}{\ell-1+\xi},
$$
ie $A_1$ contains $\frac{2D(1-\xi)}{\ell-1+\xi}$ vs less than
the `correct' number.}
In particular, for each $i\le q$ there
are at most $\tau^{2/3}|A_i|$ vertices which are $i$-useless. 
To estimate the number $u_{q+1}$
of $(q+1)$-useless vertices we count the number $e(A_{q+1},V(G)\setminus A_{q+1})$
of edges emanating from~$A_{q+1}$. We have
\begin{align*}
q|A_1|\delta(G) & -2\sum_{i=1}^q e(A_i)-2\binom{q}{2}|A_1|^2 \le 
e(A_{q+1},V(G)\setminus A_{q+1})\\
& \le 
u_{q+1}[(q-1)|A_1|+(1-\tau^{1/4})|A_1|]+(|A_{q+1}|-u_{q+1})q|A_1|
\end{align*}
which implies%
     \COMMENT{As $n\ge (\ell-1+\xi)|A_1|$ we get
$$
q|A_1|\left(1-\frac{1}{\ell-1+\xi}\right)(\ell-1+\xi)|A_1|-2\tau q|A_1|^2-q(q-1)|A_1|^2
\le -u_{q+1}\tau^{1/4}|A_1|+q|A_{q+1}||A_1|$$
and thus
\begin{align*}
u_{q+1}\tau^{1/4} & \le q|A_{q+1}|-(\ell-2+\xi)q|A_1|+2\tau q|A_1|+q(q-1)|A_1|\\
& \le q|A_1|+2\tau q |A_1|-q|A_1|+\tau^2 n\le 2\tau q|A_{q+1}|/\xi +\tau^2 n.
\end{align*}
(To get the 2nd line note that $|A_{q+1}|=n-q|A_1|= (\ell-1+\xi)|A_1|-q|A_1|+
some large constant depending on~D'$.)
Thus $u_{q+1}\le \tau^{2/3}|A_{q+1}|$.}
that the number $u_{q+1}$ of $(q+1)$-useless vertices is at most
$\tau^{2/3}|A_{q+1}|$. So in total, at most $\tau^{2/3} n$ vertices of $G$
are $i$-useless for some~$i\le q+1$. 

Given $j\neq i$, we call a vertex $x\in A_i$ \emph{$j$-exceptional}
if $x$ has at most $\tau^{1/3}|A_j|$ neighbours in $A_j$.
Thus every such $x$ is both $i$-useless and $i$-bad.%
     \COMMENT{This holds regardless whether $j=q+1$ or $j<q+1$.}
It will be important later that the number of vertices which are
$i$-useless for some~$i$ is much smaller than the number of neighbours
in~$A_j$ of a non-$j$-exceptional vertex.
By interchanging $i$-bad vertices with $i$-exceptional vertices
if necessary, we may assume that for each $i$ for which there
exist $i$-exceptional vertices, we don't have $i$-bad vertices.
Note that after we have interchanged vertices every non-$j$-exceptional
vertex still has at least $\tau^{1/3}|A_j|/2$ neighbours in~$A_j$.
Similarly, every non-$i$-bad vertex still has at most $2\tau^{1/3}|A_i|$
neighbours in~$A_i$ and every non-$i$-useless
vertex still has at least $(1-2\tau^{1/4})|A_j|$ neighbours in~$A_j$
for every~$j\neq i$.

For each index $i\le q$ in turn we now proceed as follows in order
to take care of the $i$-exceptional vertices.%
    \COMMENT{The bound on $e(A_i)$ below only works if $A_i$ is
not the small class, ie if $i<\ell$. Moreover, if $q\le \ell-2$ then
no vertex will be $(q+1)$-exceptional.}
Let $S_i\subseteq V(G)\setminus A_i$ denote the set of $i$-exceptional vertices
and assume that $s_i:=|S_i|>0$. We will choose a set $\S_i$ of $s_i$ disjoint
$z$-stars in $G[A_i]$ and interchange the star centres with the
$i$-exceptional vertices.%
     \COMMENT{Note that when doing this for some $j>i$ we will not produce
new $i$-exceptional vertices. Indeed, the only way this could happen
is if we move some vertex $x\in A_i$ to $A_j$ which was a $j$-exceptional
vertex and now becomes a (new) $i$-exceptional vertex. But this is not
possible since $x$ has lots of neighbours in $A_i$ as it is
$j$-exceptional.}
To show the existence of such stars, note that $\Delta(A_i)\le 2\tau^{1/3}|A_i|$
since by our assumption no vertex in $A_i$ is $i$-bad. Moreover, we
can bound the number of edges in $G[A_i]$ by
\begin{align*}
e(A_i)& \ge \frac{\delta(G)|A_i|-|G-(A_i\cup S_i)||A_i|-e(A_i,S_i)}{2}\\
& \ge \frac{1}{2}|A_i|\left[\left(1-\frac{1}{\ell-1+\xi}\right)n-
\left(n-\frac{n}{\ell-1+\xi}-s_i\right)-2s_i\tau^{1/3}+\frac{C}{2}\right]\\
& \ge \frac{1}{2}|A_i|\left(C/2+s_i/2 \right).
\end{align*}
We only have $C/2$ instead of $C$ in the second line since we have to compensate
for the fact that the size of the $A_i$'s is not exactly $n/(\ell-1+\xi)$.
Thus $G[A_i]$ has average degree at least $C/2+s_i/2\ge 2z$.
Lemma~\ref{disjointstars} now implies that $G[A_i]$
contains at least
$$\frac{(C/2+s_i/2)|A_i|}{8(z+1)\tau^{1/3}|A_i|}\ge s_i
$$
disjoint $z$-stars, as required.
We still denote the modified sets by $A_i$ and let $\S:=\bigcup_{i=1}^{q} \S_i$.

We now choose a set $\B$ of $|\S|$ disjoint copies of $B^*$ in~$G$, each containing precisely one
of the stars in~$\S$. Moreover, each such copy will have precisely $z$ vertices
in every $A_i$ with $i\le q$. To see that such copies exist, we will first show that 
$G[A_{q+1}]$ contains many disjoint copies of $B_1^*$, where $B_1^*$ denotes the subgraph of 
$B^*$ obtained by removing $q$ of the large vertex classes. For this, let $n_{q+1}:=|A_{q+1}|$.
Then
\begin{align}\label{mindegGq1}
\frac{\delta(G[A_{q+1}])}{n_{q+1}} & \ge
\frac{\delta(G)-|A_1\cup \dots \cup A_q|}{n}
\cdot\frac{n}{n_{q+1}}\nonumber\\
& \ge \left(1-\frac{1}{\ell-1+\xi}-\frac{q}{\ell-1+\xi}\right)
\frac{\ell-1+\xi}{\ell-q-1+\xi}\nonumber\\
& = 1-\frac{1}{\ell-q-1+\xi}\nonumber\\
& =1-\frac{1}{\chi_{cr}(B^*_1)}.
\end{align}
(The fact that $C\gg D'$ enables us to ignore the terms involving the
constant $D'$ when estimating~$n/n_{q+1}$.)%
\COMMENT{Indeed, we use that $n/n^*_1\ge n/n_1$. We have
$$(\ell-1+\xi)|A_1|=
(\ell-1+\xi)\left(\frac{n-\ell D}{\ell-1+\xi}+D\right)=
n-D(1-\xi)$$ and so $n_1=n-q|A_1|=(\ell-1-q+\xi)|A_1|+D(1-\xi)$.
Hence $$n/n_1\ge \frac{\ell-1+\xi}{(\ell-1-q+\xi)+D/|A_1|}\ge
\frac{\ell-1+\xi}{\ell-1-q+\xi}-
\frac{\ell-1+\xi}{\ell-1-q+\xi}\frac{D}{|A_1|}.
$$
Also, the intermediate 3rd line in above is: $\left(\frac{\ell-2+\xi-q}{\ell-1+\xi}-
2|B^*|\tau^{2/3}\right)\frac{\ell-1+\xi}{\ell-q-1+\xi}$}
Since $\tau^{1/5}\ll \chi_{cr}(B^*_1)-(\chi(B^*_1)-1)$, by repeated applications of the 
Erd\H{o}s-Stone theorem (see~(\ref{eqErdosStone})) we can find $\tau^{1/5}n_{q+1}$ disjoint copies 
of $B_1^*$ in $G[A_{q+1}]$. Since at most $\tau^{2/3}n_{q+1}$ vertices in $A_{q+1}$
are $(q+1)$-useless, we may assume that all of these copies of $B^*_1$ avoid the
$(q+1)$-useless vertices.
Moreover, all the stars $S\in \S$ are disjoint and 
it is easy to see that none of the vertices of such a star~$S$ can be
$i$-useless where $A_i$ is the set which originally contained~$S$.%
     \COMMENT{Otherwise it would also be $i$-bad (since $i\le q$) and we could have interchanged this
$i$-bad vertex with some $i$-exceptional vertex.}
The latter implies that each vertex of~$S$ is joined to at least
$(1-2\tau^{1/4})|A_j|$ vertices in $A_j$ for every $j\neq i$.
Thus in particular, each vertex of $S$ is joined to all vertices in almost all of the copies of $B_1^*$ 
selected above.
Altogether, the above shows that we can greedily choose the set $\B$ of $|\S|$ disjoint copies of $B^*$ as 
follows: for each copy, first choose a star $S \in \S$, then choose a copy of $B_1^*$ selected above
all of whose vertices are joined to all vertices in~$S$. If the centre of~$S$ was moved into~$A_{q+1}$,
we interchange it with some vertex in the copy of~$B^*_1$. Finally we choose the remaining vertices of $B^*$.
Let $A'_i$ be the subset of $A_i$ which contains all those vertices that do not lie
in a copy of $B^*$ in~$\B$. Note that
\begin{equation}\label{eqsizeAi'}
|A_i\setminus A_i'|\le |B^*||\S|\le |B^*|\tau^{2/3}n.
\end{equation}
After this process we have removed all the $i$-exceptional vertices for all $i\le q$.

The next step is to deal with the useless vertices (and thus also with the
$(q+1)$-exceptional vertices). For each such vertex~$x$
we will move~$x$ into another vertex class or/and we will
remove a copy of $B^*$ which contains~$x$. (We do the former
if $x$ lies in the set $U$ defined below.)
Let $U$ denote the set of all vertices in $A'_1\cup\dots\cup A'_q$
which had at most $(1-\tau^{1/4})|A_{q+1}|$ neighbours in $A_{q+1}$.
So in particular, each $u\in U$ is $i$-useless where $i\le q$ is the
index such that $u\in A'_i$. Thus $|U|\le \tau^{2/3}n$.
(Moreover, if $q=\ell-1$ then $U$ contains all the $(q+1)$-exceptional
vertices. If $q<\ell-1$ then there are no $(q+1)$-exceptional vertices.) 
Note that each $u\in U\cap A'_i$ must still have
at least $\tau^{1/3}|A'_i|$ neighbours in its own class $A'_i$.
Moreover, as in~(\ref{mindegGq1}) one can show that each $u\in U$ satisfies
\begin{equation}\label{eqNu}
|N(u)\cap A'_{q+1}|\ge \delta(G)-|A_1\cup\dots\cup A_q|-|A_{q+1}\setminus A'_{q+1}|
\stackrel{(\ref{eqsizeAi'})}{\ge} \left(1-\frac{1}{\chi_{cr}(B^*_1)}-\tau^{3/5}\right)|A'_{q+1}|.
\end{equation}
Let $A''_1,\dots,A''_{q+1}$ denote the sets obtained from the $A'_i$
by moving all the vertices in $U$ to $A'_{q+1}$.
Then~(\ref{mindegGq1}) and~(\ref{eqNu}) together with the fact that
$\tau^{1/5}\ll \chi_{cr}(B^*_1)-(\chi(B^*_1)-1)$ imply that
\begin{equation}\label{eqmindegA''q+1}
\delta(G[A''_{q+1}])\ge \left(1-\frac{1}{\chi_{cr}(B^*_1)}-\tau^{1/2}\right)|A''_{q+1}|
\ge \left(1-\frac{1}{\ell-q-1}+\tau^{1/5}\right)|A''_{q+1}|.
\end{equation}
(If $q=\ell-1$ then we will only use the first inequality in~(\ref{eqmindegA''q+1}).)
Consider the graph $K$ obtained from the complete $(q+1)$-partite
graph with vertex classes $A''_1,\dots,A''_{q+1}$ by making $A''_{q+1}$
into a clique. Let $K''$ denote the subgraph of $K$ obtained by deleting $D'$
vertices from each of the first $q$ classes and $(\ell-q)D'$ vertices from~$A''_{q+1}$.
An application of Lemmas~\ref{completeapprox} and
~\ref{completeconst} shows that by taking out at most $\ell^3|U|+D'\le 2\ell^3\tau^{2/3} n$ disjoint
copies%
     \COMMENT{We take out at most $\ell^3|U|$ disjoint copies of $B^*$ in
Lemma~\ref{completeapprox} and thus at most $\ell^2|U|$ disjoint copies of $H$.} 
of $H$ from $K''$ one can obtain a subgraph $K'''$ whose
vertex classes $A'''_1,\dots,A'''_{q+1}$ satisfy
$|A'''_i|=z|K'''|/|B^*|$ for all $i\le q$.
Moreover, each of these copies of $H$ meets~$A''_{q+1}$ in
an $(\ell-q)$-partite graph. Together with~(\ref{eqmindegA''q+1})
and the Erd\H{o}s-Stone theorem
this shows that for each of these copies of $H$ in $K''$ we can take out a
copy of $H$ from $G$ which intersects the $q+1$ vertex classes in exactly the same
way and avoids all the useless vertices. We add all these copies of $H$ in $G$
to the set~$\B$.
Adding the $\ell D'$ vertices set aside (when going from $K$ to $K''$)
to the vertex classes again we thus obtain vertex sets $A^\diamond_1,\dots,A^\diamond_{q+1}$
such that $|A^\diamond_i|=(n^\diamond-\ell D')z/|B^*|+D'$ for all $i\le q$,
where $n^\diamond:=|A^\diamond_1\cup\dots\cup A^\diamond_{q+1}|$.

By the bound in~(\ref{eqsizeAi'}) and the previous paragraph we have removed at most
$3\ell^3|B^*|\tau^{2/3}n\ll \tau^{1/3}n$ vertices so far.
Thus for all $i\le q+1$ every vertex in $A^\diamond_i$
still has at least $\tau^{1/3}|A^\diamond_j|/3$ neighbours in each other~$A^\diamond_j$
with $j\le q$ (since it is non-$j$-exceptional).
Moreover, since we have moved the vertices in $U$, for all $i\le q$
every vertex in $A^\diamond_i$
still has at least $(1-3\tau^{1/4})|A^\diamond_{q+1}|$ neighbours in~$A^\diamond_{q+1}$.
Also, $d(A^\diamond_i)\le \tau^{2/5}/2$ for all $i\le q$
(note that exchanging $i$-exceptional vertices with $i$-bad  vertices
does not affect the density too much).%
      \COMMENT{Here we have to be careful since, apart from deleting a small
number of vertices from $A_i$ we also added some new vertices when
we interchanged $i$-bad vertices with $i$-exceptional vertices
or $j$-exceptional vs (in $A_i$) with $j$-bad vertices. In the latter step
we might have added up to $\tau^{2/3}|A_i|$ vertices to $A_i$ which see
everything in $A_i$. Hence the increase in the density.}
For all $i\le q$ in turn, we now add further copies
of~$B^*$ to $\B$ in order to cover all those $i$-useless vertices which still
lie in $A^\diamond_i$. Let $U'$ denote the set of all these vertices.
Again, each such copy of $B^*$ will meet every $A^\diamond_i$ with $i\le q$ in precisely $z$ vertices.
It is easy to see that these copies of $B^*$ can be found greedily.%
     \COMMENT{Indeed, at least $\tau^{1/3}|A_j|/4-\tau^{2/3}n\gg \tau^{1/3}|A\diamond_j|/8$
of the neighbours of an $i$-useless vertex in some
$A\diamond_j$ are typical, ie not $j$-useless. So the number of these neighbours is much larger
than the number of vertices which we have already removed from $A_j$ in order to deal
with other $j'$-useless vertices, so we can proceed greedily.}
This follows similarly as before since $|U'|$ is much smaller than the number of
neighbours in any~$A^\diamond_j$ of such
a (non-$j$-exceptional) vertex $u\in U'$ and since each $u\in U'$ is joined to almost
all vertices in $A^\diamond_{q+1}$. More precisely, given $u\in U'\cap A^\diamond_k$,
let $i\le q$ with $i\neq k$ be such that $|N(u)\cap A^\diamond_i|$ is minimal.
Note that this implies that
$|N(u)\cap A^\diamond_j|\ge |A^\diamond_j|/3$ for all $j\le q$ with $j\neq i,k$.
We choose the copy of~$B^*$ containing~$u$ by first picking $z$ neighbours
of $u$ in $A^\diamond_i$ which are not $i$-useless
(this can be done since $u$ has at least $\tau^{1/3}|A^\diamond_i|/3$
neighbours in~$A^\diamond_i$). Then we pick a copy of~$B^*_1$ in~$A^\diamond_{q+1}$
which is joined to all the $z+1$ vertices chosen before 
and which also avoids all the useless vertices. Finally, we pick the remaining vertices.

Call a vertex $u\in A^\diamond_{q+1}$ \emph{worthless} if $u$ has at most
$(1-3\tau^{1/4})|A^\diamond_i|$ neighbours in~$A^\diamond_i$ for some $i\le q$.
Thus every worthless vertex is either $(q+1)$-useless
or lies in~$U$. In particular, at most $\tau^{2/3}n$ vertices are worthless.
For each worthless vertex~$u$ we will remove a copy of $B^*$ containing~$u$.
Since $u$ has at least $\tau^{1/3}|A^\diamond_i|/3$
neighbours in~$A^\diamond_i$ for each $i\le q$, it is easy to see that this
can be done if $q=\ell-1$. So suppose that $q<\ell-1$.

We now consider all the worthless vertices~$u$ in turn.
Again, we let $i\le q$ be such that
$|N(u)\cap A^\diamond_i|$ is minimal. Thus
$|N(u)\cap A^\diamond_j|\ge |A^\diamond_j|/3$ for all $j\le q$ with $j\neq i$.
Choose a set~$T_u$ of~$z$ neighbours of~$u$ in~$A^\diamond_i$.
Let $N_u$ denote the set of all those common neighbours of
the vertices from~$T_u$ in the set~$A^\diamond_{q+1}$ which are not worthless. Thus
\begin{equation} \label{Nusize}
|N_u|\ge (1-\tau^{1/5})|A^\diamond_{q+1}|.
\end{equation}
We will show that there are many disjoint copies of~$B^*_1$ in~$G[N_u]$ such that
all but one vertex class in each of these copies lie in the neighbourhood of~$u$.
We will call such a copy of~$B^*_1$ \emph{good for~$u$}.

Let $t:=\lceil 3/\xi\rceil$. Let~$K^*$ denote the complete $(\ell-q)$-partite graph
with $\ell-q-1$ vertex classes of size~$zt$ and one vertex class of size~$z_1t$.
Note that $\chi_{cr}(K^*)=\chi_{cr}(B^*_1)$. Thus
Theorem~\ref{thmKomlos} together with~(\ref{Nusize}) and
the first inequality in~(\ref{eqmindegA''q+1})
imply that $G[N_u]$ contains a $K^*$-packing which covers all but at most
$\tau^{1/6}|N_u|$ vertices. On the other hand, similarly as
in~(\ref{eqNu}) we have
\begin{align*}
|N(u)\cap N_u| & \stackrel{(\ref{Nusize})}{\ge} \left(1-\frac{1}{\chi_{cr}(B^*_1)}-\tau^{1/6}\right)|N_u|
\stackrel{(\ref{eqchicr})}{=}
\left(\frac{(\ell-q-2+\xi)zt}{(\ell-q-1+\xi)zt}-\tau^{1/6}\right)|N_u|\\
& \ge \left(\frac{(\ell-q-2)zt+2z}{|K^*|}+\tau^{1/7}\right)|N_u|,
\end{align*}
where the last inequality holds since $\xi zt\ge 2z+2|K^*|\tau^{1/7}$.
Thus there are many copies of $K^*$ in the $K^*$-packing such that $u$ is joined to at
least~$z$ vertices in all but at most one class. Each such copy of~$K^*$ gives a copy of~$B^*_1$
which is good for~$u$. Take such a copy of~$B^*_1$, exchange~$u$ with an appropriate
vertex, extend the new copy of~$B^*_1$ to a copy of~$B^*$ (which will
meet~$A^\diamond_i$ precisely in~$T_u$) and remove it.
Since there is room to spare in the calculations above, we can do this for every
worthless vertex~$u$ in turn.

Let $A^*_i$ denote the subset of all those vertices in $A^\diamond_i$ which are not covered
by some copy of $B^*$ or~$H$ in~$\B$. Then the sets $A^*_1,\dots,A^*_{q+1}$ are as required
in the lemma.%
     \COMMENT{To check that the size of $A^*_1$ is as required in (ii) let $r$
denote the number of copies of $B^*$ taken out. Thus $|A^*_1|=|A_1|-rz$
and $n^*=n-rz(\ell-1+\xi)$ and so
\begin{align*}
|A^*_1| & =\frac{z(n-\ell D')}{|B^*|}+D'-rz=
\frac{z(n^*+rz(\ell-1+\xi)-\ell D')}{|B^*|}+D'-rz\\
& =\frac{z(n^*-\ell D')}{|B^*|}+D'+\frac{z^2r(\ell-1+\xi)}{z(\ell-1+\xi)}-zr
=\frac{z(n^*-\ell D')}{|B^*|}+D'.
\end{align*}}
\endproof

We first deal with the case where $G$ looks very much like the complete $\ell$-partite
graph whose vertex class sizes are a multiple of those of the bottlegraph $B^*$ of $H$.

\begin{lemma}\label{extremal1}
Suppose that $H$ is a graph of chromatic number $\ell\ge 2$ such that
$\hcf(H)=1$. Let $B^*$ denote the bottlegraph assigned to~$H$.
Let $\xi$, $z$ and $z_1$ be as defined in~$(\ref{eqdefxi})$ and 
let $0<\tau \ll \xi, 1-\xi, 1/|B^*|$.
Let $|B^*|\ll D\ll C$ be integers such that $D$ is divisible by $2|B^*|$. Let $G$ be
a graph whose order $n\gg C,1/\tau$ is divisible by $|B^*|$ and which satisfies the
following two properties:
\begin{itemize}
\item[{\rm (i)}] $\delta(G)\ge (1-\frac{1}{\chi_{cr}(H)})n+C$.
\item[{\rm (ii)}] The vertex set of $G$ can be partitioned into $A_1,\dots,A_\ell$
such that, for all $i<\ell$, we have $|A_i|=(n-\ell D)z/|B^*|+D$ and $d(A_i)\le \tau$.
\end{itemize} 
Then $G$ has a perfect $H$-packing.
\end{lemma}
\proof
Our aim is to find a subgraph of $G$ for which it is clear that we can apply the 
Blow-up lemma to find a perfect $H$-packing. 
We first apply Lemma~\ref{exceptionalvs} with $q:=\ell-1$ and $D':=D/2$ to obtain sets
$A^*_1,\dots,A^*_\ell$.
Let $G^*$ denote the subgraph of $G$ induced by the union of all these $A^*_i$.
So $G-G^*$ has a perfect $H$-packing. It is easy to see (and follows from
Lemma~\ref{completeconst} applied with $D'=D/2$)
that the complete $\ell$-partite graph with vertex classes $A^*_1,\dots,A^*_\ell$
has a perfect $H$-packing. Since in $G^*$ each vertex in $A^*_i$ has at least
$(1-\tau^{1/5})|A^*_j|$ neighbours
in each other~$A^*_j$ the bipartite subgraph of $G^*$ between every
pair $A^*_i$, $A^*_j$ of sets is $(2\tau^{1/5},1/2)$-superregular.
Hence the Blow-up lemma implies that $G^*$ has a perfect $H$-packing.
Together with all the copies of $H$ chosen so far this yields
a perfect $H$-packing in~$G$.
\endproof

Another family of graphs having large minimum degree but not containing a perfect
$H$-packing can be obtained from a complete $\ell$-partite graph whose
vertex classes are multiples of those of the bottlegraph $B^*$ as follows:
remove all edges between the smallest vertex class ($A_\ell$ say) 
and one of the others ($A_{\ell-1}$ say), remove one vertex $x$ from $A_\ell$,
delete all the edges between $x$ and $A_1$
and add $x$ to~$A_1$, add all edges within the remainder of $A_\ell$ and
add a sufficient number of edges within $A_{\ell-1}$.
The next lemma deals with the case where $G$ is similar to this family of graphs,
although slightly more dense.

\begin{lemma}\label{extremal2}
Suppose that $H$ is a graph of chromatic number $\ell\ge 2$ such that
$\hcf(H)=1$. Let $B^*$ denote the bottlegraph
assigned to~$H$. Let $\xi$, $z$ and $z_1$ be as defined in~$(\ref{eqdefxi})$ and 
let $0<\tau \ll \xi, 1-\xi, 1/|B^*|$. Then there exists an integer $s_0=s_0(\tau,H)$
such that the following holds. Let $|B^*|\ll D\ll C$ be integers such that
$D$ is divisible by~$s_0$. Let $G$ be
a graph whose order $n\gg C,1/\tau$ is divisible by $|B^*|$ and which satisfies the
following properties:
\begin{itemize}
\item[{\rm (i)}] $\delta(G)\ge (1-\frac{1}{\chi_{cr}(H)})n+C$.
\item[{\rm (ii)}] There are disjoint vertex sets $A_1,\dots,A_{\ell-2}$ in $G$
such that $|A_i|=(n-\ell D)z/|B^*|+D$ and $d(A_i)\le \tau$ for all $i\le \ell-2$.
\item[{\rm (iii)}] The graph $G_1:=G-\bigcup_{i=1}^{\ell-2} A_i$ contains
a vertex set $A$ such that $d(A,V(G_1)\setminus A)\le \tau$.
\end{itemize}
Then $G$ has a perfect $H$-packing.
\end{lemma}
\proof
Put
$$p:=\lfloor 4/\xi\rfloor.$$
Fix further constants $\eps',d',\theta,\tau_2,\dots,\tau_p$ such that
$$
0<\eps'\ll d'\ll\theta\ll\tau\ll\tau_2\ll\tau_3\ll\dots\ll \tau_p\ll \xi,1-\xi,1/|B^*|.
$$
Let $B_1^*$ be the complete bipartite graph with vertex classes of size $z_1$ and $z$
(in other words, it is the subgraph of $B^*$ induced by its two smallest vertex classes).
Let $k_1(\eps',\theta,B^*_1)$ be as defined in Lemma~\ref{nonextremal}.
Put
$$s_0:=4k_1(p!)|B^*_1||B^*|.
$$
Let $q:=\ell-2$ and $A_{q+1}:=V(G_1)$. If $\ell\ge 3$ we first apply Lemma~\ref{exceptionalvs}
with $D'=D/2$ to obtain sets $A_1^*,\dots,A_{q+1}^*$. Let $G^*$ denote the subgraph of $G$
induced by all the~$A^*_i$ and put $n^*:=|G^*|$. Thus $G^*$ was obtained from
$G$ by taking out a small number of disjoint copies of~$H$ and
$n-n^*\le \tau^{3/5}n$. We have to show that $G^*$ has a perfect $H$-packing.
Put $G^*_1:=G[A^*_{q+1}]$ and $n^*_1:=|G^*_1|$. Note that $n^*_1$ is divisible
by~$|B^*_1|$. (In the case when $\ell=2$ we put $G^*=G^*_1=G$.)
Also, it will be crucial later on that  
\begin{equation}\label{mindegG*1}
\frac{\delta(G^*_1)}{n^*_1} 
 \ge 1-\frac{1}{\chi_{cr}(B^*_1)}-\tau^{1/2}=\frac{\xi}{1+\xi}-\tau^{1/2}.
\end{equation}
This can be proved in the same way as~(\ref{mindegGq1}), the only difference is that
we have to account for the fact that $V(G)$ and $V(G^*)$ are not quite the same.
(But since $n-n^* \le \tau^{3/5}n$, we can compensate for this
by including the error term $\tau^{1/2}$ in the above.)

Ideally, we would like to choose a perfect $B_1^*$-packing in $G_1^*$ using 
Lemma~\ref{nonextremal}. If $\ell\ge 3$ we would
like to extend each copy of $B_1^*$ to a copy of $B^*$
by adding suitable vertices in $A^*_1\cup \dots\cup A^*_q$.
Inequality~(\ref{mindegG*1}) implies that $G_1^*$ has sufficiently large minimum
degree for this. However, we cannot apply Lemma~\ref{nonextremal} directly to
$G_1^*$ since the condition (ii) is not satisfied. 
So we will consider the `almost components' of $G_1^*$ instead.

Note that if $C'\subseteq V(G^*_1)$ is such that
$d(C',V(G^*_1)\setminus C')\le \tau_p$ then $|C'|\ge \delta(G^*_1)-\tau_p n^*_1
\ge \xi n^*_1/3$. 
Let $r\le p$ be maximal such that there is a partition
$C_1,\dots,C_r$ of $V(G^*_1)$ with $d(C_j,V(G^*_1)\setminus C_j)\le \tau_r$ for
all $j\le r$. We have just seen that $r\le 3/\xi<p$ and
\begin{equation}\label{eqsizeCi}
\left(\frac{\xi}{1+\xi}-\tau_r^{1/2}\right)n^*_1\le |C_j|
\le \left(\frac{1}{1+\xi}+\tau_r^{1/2}\right)n^*_1.
\end{equation}
Moreover, $r\ge 2$ since $d(A\cap V(G^*_1),V(G^*_1)\setminus A)\le \tau_2$.

Recall that the aim is to choose a perfect $B^*_1$-packing in $G^*_1[C_j]$
(for all $j\le r$) and to extend each copy of $B^*_1$ to a copy of $B^*$
by adding suitable vertices in $A^*_1\cup \dots\cup A^*_q$.
Our choice of $r$ ensures that no $G^*_1[C_j]$ is close to an extremal
graph and thus to find a perfect $B^*_1$-packing we can argue similarly as
in the non-extremal case (c.f.~Lemma~\ref{nonextremal} and
Corollary~\ref{cornonextremal}).

More precisely, we proceed as follows. The first step is to tidy up the sets
$C_j$ to ensure that \emph{every} vertex in $C_j$ has only few neighbours
in~$G^*_1-C_j$. (The latter implies that the minimum degree of each graph $G^*_1[C_j]$
is about~$\delta(G^*_1)$.) Given a set $C'\subseteq V(G^*_1)$, put $\Cbar':=V(G^*_1)\setminus C'$.
We say that a vertex $x\in C_j$ is \emph{$j$-useless} if it has
less than $\xi|C_j|/3$ neighbours in~$C_j$. By~(\ref{mindegG*1})
each $j$-useless vertex $x$ has at least $\xi|\Cjbar|/3$ neighbours in $\Cjbar$.
Since we assumed that $d(C_j,\Cjbar)\le \tau_r$,
this implies that the number of $j$-useless vertices is at most $\tau^{3/4}_r|C_j|$. 
For all $j\le r$ we remove every $j$-useless vertex
$x\in C_j$ and add $x$ to some $C_i$ which contains at least $\xi|C_i|/3$
neighbours of~$x$. We denote by $C'_j$ the sets thus obtained from the~$C_j$.
So every vertex in $C'_j$ has at least $\xi|C'_j|/4$ neighbours in~$C'_j$.
Moreover, $d(C'_j,\overline{C'_j})\le \tau^{2/3}_r$.

We say that a vertex $x\in C'_j$ \emph{$j$-bad} if
$x$ has at least $\tau^{1/6}_r |\Cprimejbar|$ neighbours in $\Cprimejbar$.
Thus there exist at most $\tau^{1/2}_r|C'_j|$ such vertices.
For each such vertex $x$ in turn we take out a copy of $B^*_1=K_{z,z_1}$ from
$G^*_1[C'_j]$ which contains~$x$.
All these copies of $B^*_1$ can be found greedily since each
vertex in $C'_j$ has at least $\xi|C'_j|/4$ neighbours in~$C'_j$
(and thus we can apply for instance the Erd\H{o}s-Stone theorem).%
    \COMMENT{Indeed, let $N$ be a set of $\xi|C'_j|/4$ neighbours of $x$ in
$C'_j$. We distinguish two cases. Firstly, suppose that at least half of the
vertices in $N$ have at least half of there neighbours in $N$. Then $N$ is
dense and thus contains a copy of $K_{z-1},z$ which together with $x$ forms
a copy of $B^*_1$. If at least half of the vertices in $N$ have at least half
of their neighbours in $C'_j\setminus N$ then the bipartite graph between
$N$ and $C'_j\setminus (N\cup\{x\})$ is dense and thus contains
a copy of $K_{z-1,z_1}$ with $z-1$ vertices in $C'_j\setminus (N\cup\{x\})$.
Again, this gives us a copy of $B^*-1$ containing~$x$. Since the number of
$i$-bad vertices is much smaller than $\xi|C'_j|/4$,
we can argue greedily.}  
We denote by $\B'_j$ the set of all copies of $B^*_1$ chosen for the
$j$-bad vertices. Let $C''_j$ be the subset obtained from $C'_j$ in this way.
Put $n'':=|C''_1\cup \dots\cup C''_r|$.

Our next aim is to take out a \emph{bounded} number of copies of $H$
to ensure that for each $j\le r$ the size of the subset thus obtained from~$C''_j$
is divisible by~$|B^*_1|$. 
Let $D':=D/(4r)$.%
     \COMMENT{Recall that $D$ is divisible by $4r|B^*|$.}
Choose integers $t_j>0$ and $a_j$ such that
$|C_j''|=|B^*_1|t_j+a_j+4D'$ where $\sum_{i=1}^r a_j=0$ and $|a_j|<|B^*_1|$.%
      \COMMENT{Note that we can do this: Since $|B_1^*|$ divides $D'$, we can find
$a'_j$ and $t_j> 0$ with $0 \le  a'_j < |B_1^*|$ so that
$|C_j''|=|B^*_1|t'_j+a'_j+4D'$  and $\sum a'_j= s|B_1^*|$, where $s\le r$.
Note that the number of nonzero $a_j'$ is at least $s$.
Now let $a_j=a'_j-|B_1^*|$ for the first $s$ of the $a'_j$ which are nonzero
(and set $t_j=t'_j+1$ for these $j$)
and $a_j=a'_j$ and $t_j=t'_j$ for the others.} 
Note that this implies $|B^*_1|\sum t_j= n''-D$. 

We now need to distinguish the cases when $\ell\ge 3$ and $\ell=2$.
Let us first consider the case when~$\ell\ge 3$.
Let $G'$ be the graph obtained from the complete $(\ell-2)$-partite graph
with vertex classes of size $D/4$ by adding $r$ complete bipartite graphs
$K_{D'+a_j,D'}$ (with $1 \le j \le r$)
and joining all the vertices of these bipartite graphs to all the vertices
of the complete $(\ell-2)$-partite graph. So $|G'|=\ell D/4$.
An $r$-fold application of Lemma~\ref{completemove} shows that $G'$
contains a perfect $H$-packing~$\cH'$.%
     \COMMENT{We can apply this lemma since $|B^*|$ (and thus $|H|$)
divides~$D'$.}
This in turn implies that
we can greedily take out $|\cH'|= \ell D/(4|H|)$ disjoint copies of $H$ from $G$ to%
      \COMMENT{Indeed, such copies of $H$ can be found greedily since each
vertex in $A^*_i$ is joined to almost all vertices in each other $A^*_i$ as
well as to almost all vertices in each $C''_j$. Similarly each vertex in
$C''_j$ is joined to almost all vertices in each set~$A^*_i$ and to a constant
fraction of the vertices in~$C''_j$. Thus to choose a copy of $H$ which intersects
the sets $A^*_i$ as well as the set $C''_j$ (say) as desired we first
choose a suitable bipartite subgraph in $C''_j$. Then the common nbd of the
vertices in this bipartite subgraph in $A^*_i$ is almost all of $A^*_i$ and thus
we can extend the bip subgraph to the desired copy of~$H$.}
achieve that the subsets $A^\diamond_i$ and $C^\diamond_j$ thus obtained from
the sets $A^*_i$ and $C''_j$ have the following sizes
(where $n^\diamond=n^*-\ell D/4$ denotes the remaining number of vertices in the graph):%
     \COMMENT{$|A^\diamond_1|=|A^*_1|-D/4=z(n^*-\ell D/2)/|B^*|+D/2-D/4=
z(n^\diamond-\ell D/4)/|B^*|+D/4$}
\begin{equation}\label{Aidiamond}
|A^\diamond_i|=|A^*_i|-D/4 \quad \mbox{ and thus } \quad |A^\diamond_i|=z(n^\diamond-\ell D/4)/|B^*|+D/4
\end{equation}
and
\begin{equation}\label{Cdiamond}
|C^\diamond_j|=|C''_j|-2D'-a_j=|B^*_1|t_j+D/(2r).
\end{equation}
Note that every vertex in $C^\diamond_j$ still has at most
$2\tau_r^{1/6}|\overline{C^\diamond_j}|$ neighbours in~$\overline{C^\diamond_j}$.
Thus
\begin{equation}\label{eqmindegCdiamond}
\delta(G^*_1[C^\diamond_j])\stackrel{(\ref{mindegG*1})}{\ge}
\left( \frac{\xi}{1+\xi}-\tau^{1/2} \right)n^*_1 -2\tau^{1/6}_r n^*_1 
\stackrel{(\ref{eqsizeCi})}{\ge} (\xi -\tau^{1/7}_r)|C^\diamond_j|
\stackrel{(\ref{eqchicr})}{>}
\left(1-\frac{1}{\chi_{cr}(B^*_1)}\right)|C^\diamond_j|
\end{equation}
and
\begin{equation}\label{eqdensityCdiamond}
d(C^\diamond_j,\Cbar^\diamond_j)\le \tau_r^{1/7}.
\end{equation}

The arguments in the case when $\ell=2$ are similar except that
now we consider the graph $G'$ consisting
of $r$ complete subgraphs of sizes $2D'+a_j$ (where
$j=1,\dots,r$). An $(r-1)$-fold application of Lemma~\ref{completebipmove}
shows that $G'$ has a perfect $H$-packing. So we can proceed similarly
as before to obtain sets $C^\diamond_j$ which
satisfy~(\ref{Cdiamond})--(\ref{eqdensityCdiamond}).

In order to show that $G^*_1[C^\diamond_j]$  has a perfect $B^*_1$-packing we wish
to apply Lemma~\ref{nonextremal} with
$\tau_{r+1}/2$ playing the role of $\tau$ to find a blown-up $B^*_1$-cover.
(We will then use the fact that $\hcf(H)=1$ and Lemmas~\ref{completeconst} and~\ref{completeapprox}
to find a suitable perfect $B^*_1$-packing of $G_1^*[C_j^\diamond]$ which is extendable
to a perfect $H$-packing of the whole graph.) 
Inequality~(\ref{eqmindegCdiamond}) shows that
$G^*_1[C^\diamond_j]$ satisfies the requirement on the minimum degree in
Lemma~\ref{nonextremal}.
We will now check that $G^*_1[C^\diamond_j]$ also satisfies the conditions (i) and (ii)
there. So suppose that $G^*_1[C^\diamond_j]$ does not satisfy~(ii) and let
$C'\subseteq C^\diamond_j$ be such that $d(C',C^\diamond_j\setminus C')\le \tau_{r+1}/2$.
Using that $|C_i\setminus C^\diamond_i|\ll \tau_{r+1}|C_i|$
for all $i\le r$ it is easy to see that the $r+1$ sets
$C'\cap C_j$, $C_j\setminus C'$, $C_i$ ($i\neq j$)
would then contradict the choice of~$r$.%
     \COMMENT{We need $C'\cap C_j$ instead of $C'$ here since $C'\subseteq C^\diamond_j
\subseteq C'_j$ but we might have $C'\not\subseteq C_j$.}
Thus $G^*_1[C^\diamond_j]$ satisfies~(ii).
Suppose next that $G^*_1[C^\diamond_j]$ does not satisfy~(i).
Let $C'\subseteq C^\diamond_j$ be such that $|C'|=z|C^\diamond_j|/|B^*_1|$
and $d(C')\le \tau_{r+1}/2$.  
Let $x\in C'$ by any vertex which has at most $\tau_{r+1} |C'|$ neighbours in~$C'$.
Then
\begin{align*}
d_{G^*_1[C^\diamond_j]}(x) \le & \tau_{r+1} |C'|+|C^\diamond_j\setminus C'|
\le \left( \tau_{r+1}+1-\frac{z}{|B_1^*|} \right)|C^\diamond_j|
\stackrel{(\ref{eqmindegCdiamond})}{<} 
\delta(G^*_1[C^\diamond_j]),
\end{align*}
a contradiction. (In the final inequality, we also used the fact that
$1-z/|B^*_1|=\xi/(1+\xi)$.)%
     \COMMENT{Actually, we don't use the final inequality of~\ref{eqmindegCdiamond}
but the one before}
Thus we can apply Lemma~\ref{nonextremal} to $G^*_1[C^\diamond_j]$ to obtain
a $B^*_1$-packing $\B^*_j$ in $G^*_1[C^\diamond_j]$ such that the graph
$G^\diamond_j:=G^*_1[C^\diamond_j]-\bigcup \B^*_j$ (which is obtained
by removing all those vertices which lie in a copy of $B^*_1$ in~$\B^*_j$)
has a blown-up $B^*_1$-cover with parameters
$2\eps',d'/2,2\theta,k:=k_1(\eps',\theta,B^*_1)$.
Note that this blown-up $B_1^*$-cover does not necessarily yield a perfect 
$B_1^*$-packing of $G^\diamond_j$ as we need not have $\hcf(B_1^*)=1$.
This will cause difficulties later on.

Recall that, when removing the $j$-bad vertices from $C'_j$, 
we have already set aside a $B^*_1$-packing~$\B'_j$.
For each $B\in \B^*_j\cup \B'_j$ we now greedily
choose $z$ vertices in each of $A^\diamond_1,\dots,A^\diamond_q$ such that all these vertices
form a copy of $B^*$ together with~$B$.%
     \COMMENT{Indeed, this can be done greedily since every vertex in $G^*_1$ sees almost
all vertices in each $A^\diamond_i$ and all these vertices in turn see almost all
vertices in each other $A^\diamond_j$ as well as in $G^*_1$ and so on. Thus the common neighbourhood
of the vs from $B$ in some $A^\diamond_i$ is almost all of $A^\diamond_i$.}
We remove these copies of $B^*$ and still denote the subsets of the $A^\diamond_i$
obtained in this way by~$A^\diamond_i$. Also, we still denote the remaining number of vertices in
$G$ by $n^\diamond$. Then it is easy to see that the second 
equation in~(\ref{Aidiamond}) still holds for all $i$.

Consider the blown-up $B^*_1$-cover of~$G^\diamond_j$.
Let $\{ X_i^j(t) \mid 1\le t \le k, \, i=1,2\}$ be a partition of $V(G^\diamond_j)$
as in the definition of a blown-up $B^*_1$-cover (Definition~\ref{defblownupcover}).
For all $j\le r$ and all $t\le k$ we would like to apply the Blow-up lemma
in order to find a $B^*_1$-packing which covers precisely the vertices
in $X_1^j(t)\cup X_2^j(t)$. To be able to do this, we need that
the complete bipartite graph with vertex classes $X_1^j(t)$ and $X_2^j(t)$
contains a perfect $B^*_1$-packing. Clearly, the latter is the case
if $|X_1^j(t)|=z|X_1^j(t)\cup X_2^j(t)|/|B^*_1|$ and
$|X_2^j(t)|= z_1|X_1^j(t)\cup X_2^j(t)|/|B^*_1|$.  
We will now show that this can be achieved by taking out a small number of
further copies of~$H$ from~$G$. (Note that this would be much simpler to achieve
if we could assume that $\hcf(B_1^*)=1$.)

Put $D'':=D/(4rk)$ and%
     \COMMENT{Here we need that $D$ is divisible by~$4rk|B^*_1|$.}
$x^{j}(t):=(|X_1^j(t)\cup X_2^j(t)|-2D'')/|B^*_1|$.%
     \COMMENT{The definition of
blown-up $B^*$-cover implies that each $x^j(t)$ is an integer. ($x^j(t)$
counts the number of copies of $B^*_1$ in the $t$th blown-up copy of $B^*_1$
which (in each class) avoid
the `additional set' of~$2D''$ vertices.)}
If $\ell\ge 3$ consider a random partition of  $A^\diamond_i$ into $kr$ sets $A^{\diamond j}_i(t)$
($j\le r$, $t\le k$) such that $|A^{\diamond j}_i(t)|=z x^j(t)+D''$. (It is straightforward to
check that these numbers sum up to exactly $|A_i^\diamond|$.)%
      \COMMENT{
\begin{align*}
\sum_j^r \sum_t^k |A_i^{\diamond j}(t)|= &
\sum_j^r \sum_t^k (z x^j(t) +D'') = \frac{z}{|B_1^*|}
\left( \sum_j |G_j^\diamond | -2D''rk \right) +D/2 \\
=& \frac{1}{1+\xi} \left( n^\diamond -(\sum_i A_i^\diamond ) -D \right) +D/2 \\
=& \frac{1}{1+\xi} \left( n^\diamond -\ell D/2+\ell D/2-
\left(\frac{z}{|B^*|}[n^\diamond -\ell D/2] +D/2 \right) (\ell-2) -D \right) +D/2 \\
=& \frac{1}{1+\xi} \left( [n^\diamond -\ell D/2]\frac{(\ell-1+\xi )-(\ell-2 )}{\ell-1+\xi}
+\ell D/2 -(\ell-2)D/2  -D \right) +D/2 \\
=& \frac{n^\diamond-\ell D/2}{\ell-1+\xi} +D/2,
\end{align*}
as required}
Consider the complete $\ell$-partite graph $G_j(t)$ with vertex classes
$A^{\diamond j}_1(t),\dots,A^{\diamond j}_{\ell-2}(t),X_1^j(t),X_2^j(t)$.
As is easily seen, the sizes of these classes satisfy the conditions
of Lemma~\ref{completeapprox} with $D''$ playing the role of $D'$ in that lemma.%
     \COMMENT{Indeed,
$$|G_j(t)|=(\ell-2)zx^j(t)+(\ell-2)D''+|B^*_1|x^j(t)+2D''=
(\ell-1+\xi)zx^j(t)+\ell D''$$
and thus $z(|G_j(t)|-\ell D'')/|B^*|=zx^j(t)$ ie
$|A^{\diamond j}_i(t)|-D''=z(|G_j(t)|-\ell D'')/|B^*|$. The size of the
sets $X^j_1(t)$ and $X^j_2(t)$ are ok by definition of blown-up $B^*_1$
cover.}
Thus we can apply Lemma~\ref{completeapprox} and then
Lemma~\ref{completeconst} to obtain a subgraph $\tilde{G}_j(t)$ of $G_j(t)$ which is obtained
by removing a few copies of~$H$ as described there. So this gives us a collection
$\tilde{\cH}_j(t)$ of at most $\theta^{1/20} |G_j(t)|$ disjoint copies of $H$ in
$G_j(t)$ such that the subsets
$\tilde{A}^j_1(t),\dots,\tilde{A}^j_{\ell-2}(t),\tilde{X}^j_1(t),\tilde{X}^j_2(t)$
obtained from $A^{\diamond j}_1(t),\dots,A^{\diamond j}_{\ell-2}(t),X^j_1(t),X^j_2(t)$
by deleting all those vertices
which lie in some copy of $H$ in $\tilde{\cH}_j(t)$ satisfy
$$
|\tilde{A}^j_i(t)|=z\tilde{n}_j/|B^*|=|\tilde{X}^j_1(t)|
$$
for all $i\le \ell-2$ and
$$
|\tilde{X}^j_2(t)|=\xi z \tilde{n}_j/|B^*|,
$$
where $\tilde{n}_j:=|\tilde{G}_j(t)|$.
(We get the bound of $\theta^{1/20}|G_j(t)|$ copies  by observing that we
can apply Lemma~\ref{completeapprox} with $a_{\ell-1} \le \theta^{1/15} |G_j(t)|$
from Definition~\ref{defblownupcover} and $a_i=0$ for $i\le \ell-2$.)

For each copy $H'\in\tilde{\cH}_j(t)$ of $H$ in turn we greedily remove a
copy of $H$ in $G^\diamond_j\subseteq G$ which intersects the sets
$A^{\diamond j}_1(t),\dots,A^{\diamond j}_{\ell-2}(t),X_1^j(t),X_2^j(t)$
in the same way as~$H'$.
(Note that we are able to do this as Lemma~\ref{exceptionalvs}(iii) and the
fact we considered a random partition of the  $A^\diamond_i$ imply that all
vertices in  $G^\diamond_j$ are adjacent to almost all vertices in each
$A^{\diamond j}_i(t)$. Similarly, all vertices in $A^{\diamond j}_i(t)$
are adjacent to almost all vertices in each other $A^{\diamond j}_{i'}(t)$.
Moreover, the pair $(X_1^j(t),X_2^j(t))$ is $(2\eps',d'/2)$-superregular.
This enables us to construct the required number of copies of $H$ in $G^\diamond_j$
if we begin the
construction of each copy with the vertices that lie in $X_1^j(t)\cup X_2^j(t)$.
Since $\theta \gg d'$ we have to be careful that we do not destroy the
superregularity of the leftover subsets of the sets $X_i^j(t)$ in this process. We can get 
around this difficulty by considering a random red-blue partition of the vertices
as described in Section~\ref{sec:applyRG} again and removing only copies of $H$
whose vertices are all blue.)
We think of
$\tilde{A}^j_1(t),\dots,\tilde{A}^j_{\ell-2}(t),\tilde{X}^j_1(t),\tilde{X}^j_2(t)$
as the subsets obtained from
$A^{\diamond j}_1(t),\dots,A^{\diamond j}_{\ell-2}(t),X_1^j(i),X_2^j(t)$
in this way. We can now apply the Blow-up lemma to find a $B^*_1$-packing
which covers precisely the vertices in $\tilde{X}^j_1(t) \cup \tilde{X}^j_2(t)$.
The union of all these $B^*_1$-packings over all $t\le k$ and all $j\le r$ 
forms a $B^*_1$-packing $\B_1$
which covers precisely the leftover vertices of the graphs $G_j^\diamond$
with $j \le r$. If $\ell=2$ then $\B_1$ together with all the copies of $B^*_1=B^*$
and $H$ chosen earlier yields a perfect $H$-packing of~$G$.

If $\ell\ge 3$ then our aim is to extend $\B_1$ to a $B^*$-packing by
adding all the remaining vertices in the sets~$A^\diamond_i$.
So for all $i\le \ell-2$ let $A'_i$ be the subset of $A^\diamond_i$ which is left over
after removing the copies of $H$ in the union (over all $j$ and $t$) of the sets
$\tilde{\cH}_j(t)$ described above.
In order to extend $\B_1$ to a $B^*$-packing which also covers the vertices
in the sets $A'_i$, we consider the following $(\ell-1)$-partite
auxiliary graph~$J$.
The vertex classes of $J$ are $A'_1,\dots,A'_{\ell-2},\B_1$. The subgraph
of $J$ induced by $A'_1,\dots,A'_{\ell-2}$ is nothing else than the $(\ell-2)$-partite
subgraph of $G$ induced by these sets. $J$ contains an edge between
$x\in A'_j$ and $B_1^*\in \B_1$ if $x$ is joined (in $G$) to all vertices of~$B_1^*$.
Using Lemma~\ref{exceptionalvs}(iii) and the fact that we deleted only
comparatively few vertices so far, it is easy to check that in each of the
$\binom{\ell-1}{2}$ bipartite subgraphs forming $J$, every vertex is
adjacent to all but a $\tau_{r+1}$-fraction of the vertices in the other class
(with room to spare).
Thus the bipartite subgraphs are all  $(2\tau_{r+1},1/2)$-superregular.
Let $B^*_2$ be the complete $(\ell-1)$-partite graph with $\ell-2$ vertex classes of
size $z$ and one vertex class of size~1. 
The Blow-up lemma implies that $J$ has a perfect $B^*_2$-packing.
This corresponds to a $B^*$-packing (and thus also an $H$-packing) in~$G$.
Together with all the copies of $H$ chosen earlier
this yields a perfect $H$-packing in~$G$.
\endproof

The final lemma in this section deals with the remaining `extremal' possibilities: 
so $G$ contains at least one large almost independent set $A$ but does not
satisfy the conditions of either of the two previous lemmas.

\begin{lemma}\label{extremal3}
Suppose that $H$ is a graph of chromatic number $\ell\ge 3$ such that $\hcf(H)=1$.
Let $B^*$ denote the bottlegraph assigned to~$H$. Let $\xi$, $z$ and $z_1$
be as defined in~$(\ref{eqdefxi})$ and let $0<\tau \ll \tau'\ll \xi, 1-\xi, 1/|B^*|$.
Then there exists an integer $s_1=s_1(\tau,\tau',H)$
such that the following holds. Let $|B^*|\ll D\ll C$ and $1\le q\le \ell-2$
be integers such that $D$ is divisible by~$s_1$. Let $G$ be
a graph whose order $n\gg C, 1/\tau$ is divisible by $|B^*|$ and which satisfies the
following properties:
\begin{itemize}
\item[{\rm (i)}] $\delta(G)\ge (1-\frac{1}{\chi_{cr}(H)})n+C$.
\item[{\rm (ii)}] There are disjoint vertex sets $A_1,\dots,A_q$ in $G$
such that $|A_i|=(n-\ell D)z/|B^*|+D$ and $d(A_i)\le \tau$ for all $i\le q$.
\item[{\rm (iii)}] $G$ does not contain disjoint vertex sets
$A'_1,\dots,A'_{q+1}$ such that $|A'_i|=(n-\ell D)z/|B^*|+D$ and
$d(A'_i)\le \tau'$ for all $i\le q+1$.
\item[{\rm (iv)}] If $q=\ell-2$, then the graph $G_1:=G-\bigcup_{i=1}^{\ell-2} A_i$
contains no vertex set $A$ so that $d(A,V(G_1) \setminus  A) \le  \tau'$.
\end{itemize}
Then $G$ has a perfect $H$-packing.
\end{lemma}
\proof
Let $\theta:=\tau^{1/2}$.
Fix further constants $\eps',d'$ such that
$$
0<\eps'\ll d'\ll\tau\ll\tau'\ll \xi,1-\xi,1/|B^*|.
$$
Let $B_1^*$ denote the complete $(\ell-q)$-partite graph with $\ell-q-1$
vertex classes of size $z$ and one vertex class of size~$z_1$. 
Let $k_1(\eps',\theta,B^*_1)$ be as defined in Lemma~\ref{nonextremal}.
Put
$$s_1:=2k_1|B^*_1||B^*|.
$$
Let $A_{q+1}:=V(G)\setminus \bigcup_{i=1}^q A_i$. As in the proof of 
Lemma~\ref{extremal2}, we first apply Lemma~\ref{exceptionalvs} with $D'=D/2$
to obtain sets $A_1^*,\dots,A_{q+1}^*$. Let $G^*$ denote the subgraph of $G$
induced by all the~$A^*_i$ and put $n^*:=|G^*|$. Thus $G^*$ was obtained from
$G$ by taking out a small number of disjoint copies of $H$ and
$n-n^*\le \tau^{3/5}n$. Put $G^*_1:=G[A^*_{q+1}]$ and $n^*_1:=|G^*_1|$.
Note that $n^*_1$ is divisible by~$|B^*_1|$.
Moreover, as in (\ref{mindegGq1}) or (\ref{mindegG*1}) one can show that
\begin{equation}\label{mindegG*11}
\delta(G^*_1)\ge \left(1-\frac{1}{\chi_{cr}(B_1^*)}-\tau^{1/2}\right)n^*_1.
\end{equation}
Similarly as in Lemma~\ref{extremal2}, to find a perfect $H$-packing in $G^*$
our aim is to choose a perfect $B^*_1$-packing in $G^*_1$
and extend each copy of $B^*_1$ to a copy of $B^*$ by adding suitable vertices
in $A^*_1\cup \dots\cup A^*_q$. Again, we wish to apply Lemma~\ref{nonextremal}
with $\tau'/2$ playing the role of~$\tau$
to $G^*_1$ in order to do this. Thus we have to check that $G^*_1$ satisfies
the conditions of Lemma~\ref{nonextremal}. Inequality~(\ref{mindegG*11})
implies that $G^*_1$ satisfies the condition on the minimum degree.
Suppose that $G^*_1$ does not satisfy condition~(i) of Lemma~\ref{nonextremal}.
So there exists a set $A'\subseteq V(G^*_1)$ such that
$|A'|=zn^*_1/|B^*_1|$ and $d(A')\le \tau'/2$. It is easy to check that
$(1-\tau^{1/2})|A_1|\le |A'|\le (1+\tau^{1/2})|A_1|$.%
     \COMMENT{Indeed,
$$\frac{n^*_1}{|B^*_1|}=\frac{n^*-q|A^*_1|}{|B^*_1|}\approx
\frac{n^*-qzn^*/|B^*|}{|B^*_1|}= \frac{n^*}{|B^*|}\frac{(|B^*|-qz)}{|B^*_1|}
=\frac{n^*}{|B^*|},
$$
where $\approx$ means equal up to a constant depending on~$D$.
Thus we have $|A'|\approx |A^*_1|$. But $|A_1|-|A^*_1|\ll \tau^{1/2}|A_1|$.}
Thus by changing a small number of vertices we obtain a set
$A''\subseteq V(G)\setminus (A_1\cup\dots\cup A_q)$ such that
 $d(A'')\le 2\tau^{1/2}+\tau'/2\le \tau'$.
But then the sets $A_1,\dots,A_q,A''$ contradict condition~(iii)
in Lemma~\ref{extremal3}. So $G^*_1$ satisfies condition~(i) of
Lemma~\ref{nonextremal}. 
Next suppose that $q=\ell-2$ and $G^*_1$ does not satisfy
condition~(ii) of Lemma~\ref{nonextremal}.
So $G_1^*$  contains a set $A$ with $d(A,V(G_1^*) \setminus A) \le \tau'/2$.
Since $V(G_1)$ and $V(G_1^*)$ are almost the same, 
this in turn implies that the graph $G_1$ defined in (iv) 
contains a set $A$ with $d(A,V(G_1) \setminus A) \le 2\tau'/3$, a contradiction
to the assumption in (iv).
Thus we may assume that $G^*_1$ also satisfies condition~(ii) of
Lemma~\ref{nonextremal}. So we can apply Lemma~\ref{nonextremal} to~$G^*_1$.
This shows that we can take out a small number of copies of $B^*_1$ to obtain
a subgraph of $G^*_1$ which has a blown-up $B^*_1$-cover. 
We can then proceed similarly as in the final part of the proof of Lemma~\ref{extremal2}.
\endproof

\section{Proof of Theorem~\ref{thmmain}}\label{sec:thmproof}

We will now combine the results of Sections~\ref{sec:nonextremal}
and~\ref{sec:extremal} to prove Theorem~\ref{thmmain}.
Fix constants
$$0<\eps'\ll d'\ll \theta\ll \tau_1\ll\dots\ll \tau_{\ell-1}\ll \xi,1-\xi,1/|B^*|.
$$
Let $D\gg |B^*|$ be an integer satisfying the conditions
in Lemmas~\ref{extremal1}--\ref{extremal3}. Let $C\gg D,1/\tau$.
By taking out at most $\ell-2$ disjoint copies of $H$ from $G$ if necessary, we may
assume that the order $n$ of our given graph $G$ is divisible by~$|B^*|$.
(The existence of such copies follows
from the Erd\H{o}s-Stone-theorem.)
We are done if $G$ satisfies conditions~(i) and~(ii) of
Corollary~\ref{cornonextremal} with $\tau_1$ playing the role of~$\tau$.

So suppose first that~$G$ violates~(i). Thus there is some set
$A\subseteq V(G)$ of size $zn/|B^*|$ such that $d(A)\le \tau_1$.
Choose $q\le \ell-1$ maximal such that there are disjoint sets
$A'_1,\dots,A'_q\subseteq V(G)$ with $|A'_i|=zn/|B^*|$ and
$d(A'_i)\le \tau_q$ for all $i\le q$. So $q\ge 1$ by our assumption.
By removing a constant number of vertices from each $A'_i$
we obtain subsets $A_i$ with $|A_i|=z(n-\ell D)/|B^*|+D$
and $d(A_i)\le 2\tau_q$. If $q=\ell-1$ then Lemma~\ref{extremal1} shows
that $G$ has a perfect $H$-packing. Thus we may assume that $q\le \ell-2$.
But then we can apply either Lemma~\ref{extremal2} with $\tau_{q+1}$ playing the 
role of $\tau$ or Lemma~\ref{extremal3}
with $2\tau_q$ playing the role of $\tau$ and $\tau_{q+1}$ playing the
role of~$\tau'$.

If $G$ satisfies condition~(i) in Corollary~\ref{cornonextremal} but
violates~(ii) then we are done by~Lemma~\ref{extremal2} applied
with $\tau:=\tau_1$. This completes the proof of Theorem~\ref{thmmain}.
\endproof

\section{Acknowledgement}

We would like to thank Oliver Cooley for his comments on an earlier version
of this manuscript.

{\footnotesize
\bigskip\obeylines\parindent=0pt
Daniela K\"uhn \& Deryk Osthus
School of Mathematics
Birmingham University
Edgbaston
Birmingham B15 2TT
UK
{\it E-mail addresses}: {\tt \{kuehn,osthus\}@maths.bham.ac.uk}
}

\end{document}